\documentclass{amsart}
\usepackage{xypic}
\usepackage{graphicx, amsmath, amssymb, amsthm}
\usepackage{tikz-cd}
\usepackage{hyperref} 
\usepackage{mathrsfs}
\usepackage{rotating}

\newcommand{\NSFThree}{NSF Grant DMS-1509652}

\newcommand{\R}{{\mathbb R}}

\newcommand{\smashover}[1]{\underset{#1}{\wedge}}
\newcommand{\timesover}[1]{\underset{#1}{\times}}

\DeclareMathOperator{\Aut}{Aut}

\DeclareMathOperator{\Map}{Map}

\DeclareMathOperator*{\colim}{colim}

\newcommand{\tr}{tr}

\newcommand{\m}[1]{{\protect\underline{#1}}}

\newcommand{\cc}[1]{\mathcal #1}
\newcommand{\cC}{\cc{C}}
\newcommand{\ccD}{\cc{D}}
\newcommand{\cF}{\cc{F}}
\newcommand{\cG}{\cc{F}}
\newcommand{\cO}{\cc{O}}

\newcommand{\cP}{\cc{P}}

\newcommand{\cE}{\cc{E}}

\newcommand{\aC}{\cC}

\newcommand{\aO}{\cc{O}}
\newcommand{\aS}{\cc{S}}

\newcommand{\aW}{\cc{W}}

\newcommand{\Sp}{\mathcal Sp}

\newcommand{\Set}{\mathcal Set}

\newcommand{\Ninfty}{N_\infty}

\newcommand{\Top}{\mathcal T\!op}

\newcommand{\Alg}{\mathcal Alg}
\newcommand{\cOrb}{\mathcal Orb}

\newcommand{\mTop}{\m{\mathcal T\!o}p}
\newcommand{\TopG}[1][G]{\Top^{#1}}

\newcommand{\SW}{\mathcal SW}
\newcommand{\FreeO}[1][\cO]{\mathbb P_{#1}}
\newcommand{\redFreeO}[1][\cO]{\tilde{\mathbb P}_{#1}}
\newcommand{\WirtO}{\omega_{\cO,H}}

\mathchardef\mhyphen=45

\newtheorem{theorem}{Theorem}[section]
\newtheorem{lemma}[theorem]{Lemma}
\newtheorem{corollary}[theorem]{Corollary}

\newtheorem{definition}[theorem]{Definition}
\newtheorem{proposition}[theorem]{Proposition}

\newtheorem{remark}[theorem]{Remark}

\newtheorem{notation}[theorem]{Notation}

\newcommand{\defemph}[1]{\textbf{#1}} 
\newcommand{\Ho}{\textrm{Ho}}
\newcommand{\op}{\textrm{op}}
\newcommand{\htp}{\simeq}
\newcommand{\bZ}{\mathbb{Z}}
\newcommand{\fib}{\textrm{fib}}

\newcommand{\cK}{\mathcal{K}}

\newcommand{\GTop}{\TopG[G]}

\newcommand{\sma}{\wedge}
\newcommand{\bN}{\mathbb{N}}
\newcommand{\bO}{\mathbb{O}}

\newcommand{\bK}{\mathbb{K}}
\newcommand{\bQ}{\mathbb{Q}}

\newcommand{\Orb}{\textrm{Orb}}
\newcommand{\Fun}{\textrm{Fun}}

\title[Incomplete equivariant stable categories]{Equivariant stable
  categories for incomplete systems of transfers}

\begin{document}

\author[A.~J.~Blumberg]{Andrew~J. Blumberg}
\address{University of Texas \\ Austin, TX 78712}
\email{blumberg@math.utexas.edu}
\thanks{A.~J.~Blumberg was supported in part by NSF grant DMS-1151577}

\author[M.~A.~Hill]{Michael~A. Hill}
\address{University of California Los Angeles\\ Los Angeles, CA 90095}
\email{mikehill@math.ucla.edu}
\thanks{M.~A.~Hill was supported in part by {\NSFThree}}

\begin{abstract}
In this paper, we construct incomplete versions of the equivariant
stable category; i.e., equivariant stabilization of the category of
$G$-spaces with respect to incomplete systems of transfers encoded by
an $\Ninfty$ operad $\aO$.  These categories are built from the
categories of $\aO$-algebras in $G$-spaces.  Using this operadic
formulation, we establish incomplete versions of the usual structural
properties of the equivariant stable category, notably the tom Dieck
splitting.  Our work is motivated in part by the examples arising from
the equivariant units and Picard space functors.
\end{abstract}

\maketitle

\tableofcontents

\section{Introduction}

The foundation of stable homotopy theory is the study of the stable
category.  Boardman's original construction of the stable category
involved a process of formal stabilization of the category of finite
CW complexes and then the addition of colimits (e.g.,
see~\cite{Vogt}).  Subsequently, the needs of calculation drove a
lengthy period of focus on explicit models of the stable category as
the homotopy category of various constructions of categories of
spectra or (for the connective subcategory) as grouplike algebras in
spaces over the $E_\infty$-operad.  In the equivariant setting, until
very recently essentially the only constructions of the equivariant
stable category were given in terms of homotopical categories of
equivariant spectra (e.g., see~\cite{LMS} and~\cite{MandellMay}).

Recently, there has been renewed interest in formal descriptions of
the equivariant stable category in terms of structural properties.  A
substantial amount of this activity has been motivated by diagrammatic
characterizations of the stable category that are suggested by the
framework of $\infty$-categories; for example, regarding the
equivariant spectral category as spectral presheaves on the Burnside
category (e.g., see~\cite{GuillouMaypresheaves}
and~\cite{BarwickClan}) or suitable diagrams of geometric fixed-points
glued via the Tate construction (e.g., see~\cite{Glasman}
and~\cite{AbramKriz}).  A significant aspect of establishing the
correctness of these descriptions is grappling with the formal
properties of the equivariant stable category, which are considerably
more subtle than those of the stable category, especially when $G$ is
a compact Lie group.

We take the perspective that the equivariant stable category should be
thought of as characterized by equivalently
\begin{enumerate}
\item the existence of transfer maps, 
\item the Wirthmuller isomorphisms, 
\item or the tom Dieck splitting of suspension spectra. 
\end{enumerate}
For instance, as part of an effort to understand equivariant infinite
loop space theory for infinite compact Lie groups, the first author
gave a characterization of the equivariant stable category in terms of
the existence of Wirthmuller isomorphisms~\cite{Blumthesis}.

These efforts to understand the structural properties of the
equivariant stable category received a boost from the foundational
work done as part of the Hill--Hopkins--Ravenel solution to the Kervaire
invariant one problem~\cite{HHR}.  Notably, the introduction of the
multiplicative norm maps provided a new unifying language to think
about the formal structure of the equivariant stable category.
Previously, the authors studied the question of the way in which
operadic algebras in spectra and spaces for various kinds of
equivariant $E_\infty$ operads admit ``norm maps''; in spectra, these
norm maps are the multiplicative norms, and in spaces, they are
transfers~\cite{BHNinfty}.  The result of that work was to introduce
intermediate versions of equivariant commutative ring structures,
structured by $\Ninfty$ operads, interpolating between algebras over
the naive $E_\infty$ operad regarded as $G$-trivial and the
``genuine'' $E_\infty$ $G$-operad that describes classical equivariant 
commutative ring spectra.  In subsequent work, the authors gave an
algebraic characterization of compatible systems of norms or
transfers, called indexing systems~\cite{BHIncomplete}.

The purpose of this paper is to construct and study the equivariant
stable categories that are determined by the indexing systems specified
by arbitrary $\Ninfty$ operads $\cO$.  Of course, when the indexing
system comes from an equivariant little disks operad, there already
exist constructions of incomplete equivariant stable categories
described in terms of categories of orthogonal $G$-spectra specified
by an incomplete universe $U$.  However, a surprising result
of~\cite{BHNinfty} is that not all $\Ninfty$ operads are equivalent to
little disks operads; further examples of this phenomenon are given
in~\cite{Rubin2}.  Therefore, the work of this paper constructs new
kinds of equivariant stable categories that do not seem to admit
models using existing technology.  Our approach is to build a category
of spectra completely removed from the universes, considering only the
combinatorial data encoded by the \(\Ninfty\) operads. 
The way we do this has the interesting consequence that the structural
properties of the stable category are directly exposed.

The starting point for our description of the $\cO$-stable category is
the category $\GTop_\ast[\redFreeO]$ of $\cO$-algebras in based
$G$-spaces, where here $\redFreeO$ denotes the monad on based spaces 
associated to $\cO$.  We will generically refer to an object of
$\GTop_\ast[\redFreeO]$ as an $\Ninfty$ space.  Recall that a
$\cO$-algebra $B$ is grouplike if $\pi_0^H(B)$ is a group for every
closed subgroup $H \subset G$.  A map $A \to B$ of $\cO$-spaces is a
group completion if $B$ is grouplike and each induced map of
fixed-point spaces $A^H \to B^H$ is a (non-equivariant) group
completion.

We take the point of view that grouplike $\Ninfty$-spaces model the
connective part of the $\cO$-stable category, and we will produce the
rest of the category by formal stabilization.  To carry this out, we
begin by constructing a model structure on the connective objects.
The following theorem is standard.

\begin{theorem}
Let $\cO$ be an $\Ninfty$ operad.  There is a cofibrantly generated
model structure, the {\em standard model structure}, on
$\GTop_\ast[\redFreeO]$ where the weak equivalences and fibrations are
determined by the forgetful functor $\GTop_\ast[\redFreeO] \to
\GTop_\ast$.
\end{theorem}

The first basic observation is that we can group-complete an
$\cO$-algebra via the nonequivariant delooping.  Specifically, we have
the following result:

\begin{proposition}
Let $A$ be a cofibrant $\cO$-algebra.  Then the natural map 
\[
A \to \Omega (S^1 \otimes A)
\]
is a group completion, where $S^1 \otimes (-)$ denotes the based
tensor and $\Omega(-)$ the based cotensor with $S^1$.
\end{proposition}

Here we are using the fact that the based space $S^1$ is a model of
the bar construction.  Following~\cite{SagaveSchlichtkrull}, we now
localize this model structure:

\begin{definition}
Let $\cO$ be a reduced $\Ninfty$ operad.  The {\em $\cO$-stable model
  structure} on the category $\GTop_\ast[\redFreeO]$ has weak equivalences the
maps $A \to B$ such that
\[
\Omega (S^1 \otimes A) \to \Omega (S^1 \otimes B)
\]
is a weak equivalence.  The fibrant objects are the grouplike
$\cO$-algebras.
\end{definition}

We refer to $\Ho(\GTop_\ast[\redFreeO])$, formed with respect to the
$\cO$-stable model structure, as the {\em connective $\cO$-stable
  category}.  Note that the transfers that arise from the $\Ninfty$
operad impose structure on the homotopy groups of an algebra.

\begin{proposition}
Let $\cO$ be a reduced $\Ninfty$ operad and let $B$ be a grouplike
$\cO$-algebra.  Then the homotopy groups $\{\pi_n^H(B)\}$ naturally
extend to an $\cO$-Mackey functor. 
\end{proposition}

We now turn to the formal addition of nonconnective objects to get a
model for $\cO$-spectra.  For any category $\aC$ tensored and
cotensored over based spaces, we can define the category $\Sp(\aC)$ of
spectra in $\aC$ as the collections of objects $\{B_i\}$ equipped with
structure maps $S^1 \otimes B_n \to B_{n+1}$.  This definition is
motivated by the fact that when $\aC$ is a model category, the based
tensor of a cofibrant object with $S^1$ models the Quillen suspension
on $\Ho(\aC)$.

\begin{definition}
Let $\cO$ be a reduced $\Ninfty$ operad.  The category of
$\cO$-spectra \(\Sp^{G}_{\cO}\) is the category $\Sp^{G}_{\cO}$ of
spectrum objects in $\GTop_\ast[\redFreeO]$. 
\end{definition}

Standard arguments show that this has a stable model structure lifted
from the standard model structure on $\GTop_\ast[\redFreeO]$,
where the {\em stable equivalences} are detected by maps into
$\Omega$-spectra or by stable homotopy groups.

\begin{theorem}
The category of $\cO$-spectra has a cofibrantly generated
model structure where the weak equivalences are detected by the stable
equivalences.  We refer to $\Ho(\Sp^{G}_{\cO})$ as the {\em
  $\cO$-stable category}.
\end{theorem}

The functor $\Omega (S^1 \otimes (-))$ induces a fully-faithful
embedding of the connective $\cO$-stable category into the
$\cO$-stable category.  In addition, there is an adjoint pair of
functors  
\[
\xymatrix{
\Sigma^{\infty}_{\cO} \colon \GTop_{\ast} \ar@<1ex>[r] & \ar@<1ex>[l] \Sp^{G}_{\cO}\colon \Omega^{\infty}_{\cO},
}
\]
where $\Sigma^{\infty}_{\cO}$ is the composite of the free
$\cO$-algebra functor and the suspension spectrum.

Notice that this is conceptually a different approach than the
standard construction of equivariant spectra; we do not invert
representation spheres, but just the trivial spheres in operadic
algebras.  The fact that this works is precisely because grouplike
$\Ninfty$-algebras in spaces are already stable objects, and we are
simply formally adding shifts.  As a consequence, an advantage of this
approach to the construction of the $\cO$-stable category is that
basic results follow directly from analysis of the free $\cO$-algebra
in spaces.  For example, we obtain the following incomplete version of
the tom Dieck splitting by studying the fixed points of the free
$\cO$-algebra.

\begin{theorem}
Let $X$ be a $G$-CW complex. Then we have a natural equivalence 
\[
\big(\Sigma^{\infty}_{\cO} X_+\big)^G\simeq
\bigvee_{G/H\in\pi_0\cO(G)} 
EW_G(H)_{+}\smashover{W_G(H)} \Sigma^{\infty}_+X^H. 
\]
(Here $G/H \in \pi_0 \cO(G)$ denotes the isomorphism classes of
admissible $G$-sets for the indexing system determined by $\cO$.)
\end{theorem}

Equivalently, we also have an incomplete version of the Wirthm\"uller
isomorphism.

\begin{theorem}
Let $G/H$ be an admissible set for $\cO$.  Then the functor $G_+
\sma_H (-)$ is both the left and the right adjoint to the forgetful
functor $\iota^*_H$ on the $\cO$-stable category.
\end{theorem}

These results provide an axiomatic characterization of the
$\cO$-stable category.  
We also have the following consistency check.  For the
next result, recall that a map $\cO \to \cO'$ of $\Ninfty$ operads is
a weak equivalence when each constitutent map $\cO(n) \to \cO'(n)$ is
an equivalence of $G \times \Sigma_n$-spaces.

\begin{proposition}
Let $\cO$ be weakly equivalent to a little disk operad $\ccD(U)$ for a
$G$-universe $U$.  Then the underlying $\infty$-category for the
stable model structure on $\Sp^{G}_{\cO}$ is equivalent to the
underlying $\infty$-category for the stable model category of
orthogonal $G$-spectra on $U$.  The homotopy category of orthogonal
$G$-spectra for the universe $U$ is equivalent to the $\cO$-stable
category.
\end{proposition}

To explain the restriction to little disks, it is illuminating to
indicate why a simple approach using \(\Gamma_{G}\)-spaces does not in
general suffice to describe the $\cO$-stable category.  In order to
work with $\Gamma_G$ spaces, we can consider analogues of Shimakawa's
special condition, restricting attention only to the maps coming from
the admissible \(G\)-sets.  That is, we would require the maps
\[
X(T) \to F(T, X(1^+))
\]
to be $G$-equivalences for admissible $G$-sets $T$.  When all of the
admissible sets for $\cO$ are generated by the admissible sets for
\(G\) itself, this captures the correct homotopy theory.  However, for
certain other \(\Ninfty\) operads (e.g., the linear isometries operad
considered below), all of the interesting structure in the associated
indexing system takes place at proper subgroups of $G$.

Finally, we turn to explain in more detail our motivation for the work
of this paper.  As we mentioned above, our interest in $\cO$-stable
categories arose from the surprising fact that there exist exotic
$\cO$-stable categories that do not arise from categories of
equivariant orthogonal spectra for any universe.  In fact, there are
interesting examples of this phenomenon that are simple enough to
describe completely explicitly.

Consider the \(C_{4}\)-universe \(U=\infty(\mathbb
R\oplus\lambda)\), where \(\lambda\) is the unique irreducible
\(2\)-dimensional real representation. We can consider
\(C_{4}\)-spectra indexed on \(U\), and the homotopy coefficient systems
here naturally have transfer maps from the trivial group to \(C_{2}\)
and to \(C_{4}\), but no transfer from \(C_{2}\) to \(C_{4}\). The
category of \(C_{4}\)-spectra indexed on \(U\) is a symmetric monoidal
category, and commutative monoids are equivalent to operadic algebras
for the linear isometries operad structured on \(U\). A very curious
fact which motivated the development of \(\Ninfty\) operads is that
the little disks and linear isometries operads for \(U\) are distinct
and moreover, the linear isometries operad for \(U\) is not equivalent
to the little disks operad for any \(C_{4}\) universe. If we try to
form \(GL_{1}\) or related constructions, then we still produce a
group-like operadic algebra, but the only way to build an associated
spectrum is to restrict the structure to little disks on a trivial
\(C_{4}\)-universe, throwing away structure.  

This example arises from geometric considerations.  Specifically,
Manolescu's \(Pin(2)\)-equivariant approach to the \(11/8\)-conjecture
considered \(Pin(2)\)-equivariant stable homotopy indexed on a
particular incomplete universe \cite{ManolescuK} (see also
\cite{HLSX}). When restricted to \(C_{4}\subset S^{1}\subset Pin(2)\),
this universe is exactly the universe \(U\). 

We should note that despite the obvious connections to equivariant
infinite loop space theory, we do not construct a delooping machine
for arbitrary $\Ninfty$ spaces, although we think this is a very
interesting topic of study.  Also, a distinct disadvantage of the
formalism we study is that the multiplicative structure is less
explicit.  There are models for the $\cO$-stable category based on
generalizations of the excisive functor~\cite{Blumthesis} and equivariant
symmetric spectra~\cite{Mandellsymm, Hausmannsymm} approaches.  We
intend to return to these kinds of approaches in future work with
Dotto and Mandell.

\subsection*{Acknowledgments.}  
The authors would like to thank Emanuele Dotto, Mike Hopkins, Mike
Mandell, Peter May, Jonathan Rubin, and Andrew Smith for many helpful
conversations.  This project was made possible by the hospitality of
AIM and the Hausdorff Research Institute for Mathematics at the
University of Bonn.

\section{$\Ninfty$-operads and algebras}

In this section, we briefly review the notion of an $\Ninfty$ operad
and record some structural properties of categories of algebras over
these operads.  We refer the reader to~\cite{BHNinfty, BHIncomplete}
for more detailed discussion.  Throughout the following discussion, we
as usual work in the category of compactly-generated weak Hausdorff
spaces.  

\begin{notation}
We write \(\TopG\) for the category of \(G\)-spaces and
equivariant maps between them, viewed as a category enriched in
spaces.  We write \(\mTop\) for the category of \(G\)-spaces and all
maps between them, viewed as a category enriched in \(\TopG\). In both cases, we denote the corresponding pointed categories with a subscript of an asterisk.
\end{notation}

\subsection{$\Ninfty$-operads and indexing systems}

Equivariant stable categories are controlled by certain kinds of
$G$-operads.  We will assume that the $G$-operads we consider have
spaces which are of the homotopy type of $G$-CW complexes.

\begin{definition}
An $\Ninfty$ operad is an operad $\aO$ in $G$-spaces such that
\begin{enumerate}
\item The space $\aO(0)$ is $G$-contractible,
\item the action of $\Sigma_n$ on $\aO(n)$ is free, and
\item the space $\aO(n)$ is a universal space for a family of
  subgroups of $G \times \Sigma_n$ which contains all subgroups of the 
  form $H \times \{1\}$ for $H \subseteq G$.
\end{enumerate}
We say that an $\Ninfty$-operad is {\em reduced} if $\aO(0)$ is a
point.
\end{definition}

A central result about $\Ninfty$ operads is that the homotopy category
of $\Ninfty$ operads is essentially algebraic.

\begin{definition}
A weak equivalence of $\Ninfty$ operads is a map $\aO \to \aO'$ that
induces weak equivalences $\aO(n) \to \aO'(n)$ of $G \times
\Sigma_n$-spaces for each $n$.
\end{definition}

Inverting these weak equivalences, the homotopy category of $\Ninfty$
operads is equivalent to a poset~\cite{BHNinfty, GutWhite, Rubin,
  Pereira} described in terms of {\em indexing systems}. 

\begin{definition}
A {\em symmetric monoidal coefficient system} is a contravariant
functor from the orbit category of $G$ to the category of small
symmetric monoidal categories and strong monoidal functors.
\end{definition}

The canonical example of a symmetric monoidal coefficient system is
$\m{\Set}$, which assigns to $G/H$ the category of finite
$H$-sets.

\begin{definition}
An {\em indexing system} $C$ is a symmetric monoidal full sub-coefficient
system of $\m{\Set}$ that contains all trivial sets and is
closed under 
\begin{enumerate}
\item passage to finite limits and
\item self-induction, in the sense that if $H/K \in C(H)$ and $T \in
  C(K)$, then $H \timesover{K} T \in C(H)$.
\end{enumerate} 
\end{definition}

Indexing systems form a poset ordered by inclusion, and the theorem is
that the homotopy category of $\Ninfty$ operads is equivalent to the
poset of indexing systems.  Another interesting characterization of
indexing systems and hence $\Ninfty$ operads uses the theory of
polynomials:

\begin{theorem}
There is an isomorphism of categories between the poset of indexing
systems and the poset of wide, pullback stable, finite coproduct
complete subcategories of the category of finite $G$-sets.
\end{theorem}

At the extremes, there are two distinguished $\Ninfty$ operads,
corresponding to the initial and terminal indexing systems:
$\cO^{tr}$, the homotopy initial $\Ninfty$ operad, which has the
homotopy type of an ordinary, non-equivariant $E_{\infty}$-operad
viewed as a $G$-operad with trivial $G$-action, and $\cO^{gen}$, the
homotopy terminal $\Ninfty$ operad, which is a $G$-$E_{\infty}$-operad
in the sense of~\cite{LMS}.  Since the $\Ninfty$ operad $\cO^{tr}$ is
initial in the homotopy category of $\Ninfty$ operads, for any $\cO$
we have a canonical map in the homotopy category
\[
\cO^{tr}\to\cO,
\]
which can be represented on the point-set level by the inclusion
$\cO^{G} \to \cO$.

\begin{remark}
Any $\Ninfty$ operad is weakly equivalent to one which is reduced,
and hence the homotopy category of $\Ninfty$ operads can be
completely described in terms of reduced operads.  In particular, the
$\Ninfty$ Barratt-Eccles operads described in~\cite{Rubin} are reduced.
Since it is technically convenient, we will assume throughout the
remainder of the paper that our operads are reduced.
\end{remark}

\subsection{Homotopical categories of $\aO$-algebras}
\subsubsection{Point-set categories}

For an $\Ninfty$ operad $\aO$, we are interested in the category of
$\aO$-algebras in $G$-spaces and based $G$-spaces. 

\begin{definition}
Let $\FreeO$ denote the monad on $\GTop$ associated to the operad $\cO$;
to be precise, this is specified as usual by the formula
\[
\FreeO X = \coprod_{i \geq 0} \aO(i) \timesover{\Sigma_{i}} X^{\times i}.
\]
\end{definition}

The reduced monad $\redFreeO$ on $\GTop_\ast$ is specified by defining
$\redFreeO X$ to be the quotient of $\FreeO X$ by basepoint
identifications that ensure that the basepoint in $X$ coincides with
the basepoint induced by $\cO(0)$.  These identifications are
specified as follows: there are maps $\sigma_i \colon \cO(n) \to
\cO(n-1)$ that evaluate on the identity everywhere except at $i$ and
maps $s_i \colon X^n \to X^{n+1}$ that insert $\ast$ at position $i$.
We impose the equivalence relation $(o, s_i x) \sim (\sigma_i o, x)$
for $o \in \cO(n)$ and $x \in X^{n-1}$.  (See for example~\cite[\S
  4]{MayRant}.)

Note that the categories of monadic algebras $\GTop[\FreeO]$ and
$\GTop_\ast[\redFreeO]$ are equivalent; this follows from the
isomorphism of monads $\redFreeO \circ (-)_+ \cong \FreeO$ on unbased spaces
by~\cite[II.6.1]{EKMM}.  (See also~\cite[4.4]{MayRant} for further
discussion of this comparison.)

The usual arguments (e.g., see~\cite[I.7.2]{EKMM}) show that
$\FreeO$ preserves reflexive coequalizers, and the category
$\GTop_{\ast}[\redFreeO]$ is complete and cocomplete.  Limits are
created in $\GTop_{\ast}$, as are some colimits.

\begin{proposition}\label{prop:georeal}
Sifted colimits in $\GTop_{\ast}[\redFreeO]$ are created in
$\GTop_{\ast}$.  
\end{proposition}

There is an enrichment of the category $\GTop_\ast[\redFreeO]$ in
based $G$-spaces, where we consider all maps, not just equivariant ones.  

\begin{definition}
If \(A\) and \(B\) are \(\cO\)-algebras, then let
\[
\mTop_\ast[\redFreeO](A,B)
\]
be the \(G\)-space of all (not necessarily equivariant) maps
\[
A\to B
\]
which commute with monadic structure maps. Let
\(\mTop_\ast[\redFreeO]\) be the category enriched over
\(\Top^{G}_\ast\) with objects \(\cO\)-algebras and these hom
\(G\)-spaces.
\end{definition}

The two enrichments are connected by passage to $G$-fixed points.
Using the equivariant enrichment, we can talk about the standard
tensoring and cotensoring over \(G\)-spaces and based $G$-spaces.

\begin{proposition}
The category $\GTop_\ast[\redFreeO]$ is tensored and cotensored over
based \(G\)-spaces. These are defined via the adjunctions  
\[
\mTop_\ast[\redFreeO]\big(X\otimes A,B\big)\cong
\GTop_\ast \big(X,\mTop_\ast[\redFreeO](A,B)\big)\cong \mTop_\ast[\redFreeO]\big(A,F(X,B)\big). 
\]
\end{proposition}

Finally, on free algebras, the tensoring operation is especially
simple, while the underlying based \(G\)-space for the cotensoring is
especially simple.

\begin{proposition}
Let $X$ be a based $G$-space and \(A=\redFreeO(Y)\).  Then
\[
X\otimes A\cong \redFreeO(X\wedge Y).
\] 
For any object $B$ in $\GTop_\ast[\redFreeO]$, the underlying
\(G\)-space for \(F(X,B)\) is the ordinary \(G\)-space of maps
\(\mTop_\ast(X,B)\).
\end{proposition}

\subsubsection{Model structures}

We will use several model structures on \(\GTop_\ast[\redFreeO]\),
beginning with one which arises quite naturally from work of
Schwede--Shipley~\cite{SchwedeShipley}.  Specifically, since all
objects in $\GTop$ are fibrant, the following theorem follows easily.

\begin{theorem}\label{thm:ModelStructure}
There is a compactly generated $G$-topological model structure on
\(\GTop_\ast[\redFreeO]\) in which the fibrations and weak equivalences are
determined by the forgetful functor to $\GTop_\ast$.  The cofibrant
objects are retracts of filtered colimits of pushouts attaching free
cells.  We call this the {\em standard model structure}.
\end{theorem}

The assertion that the model structure is $G$-topological means that the
enrichment and (co)tensor structure of $\GTop_\ast[\redFreeO]$ satisfies the
analogue of Quillen's SM7.

For a map of $\Ninfty$ operads $\iota \colon \aO \to \aO'$, there are
associated maps of monads $\FreeO \to \FreeO[\cO']$ and $\redFreeO \to
\redFreeO[\cO']$.  We denote by 
$(\iota_{\redFreeO}^{\redFreeO[\cO']},
\iota_{\redFreeO[\cO']}^{\redFreeO})$ the adjunction between
$\GTop_\ast[\redFreeO[\cO]]$ and $\GTop_\ast[\redFreeO[\bO']]$.  Since
the right adjoint evidently preserves fibrations and weak
equivalences, the following lemma is immediate.

\begin{lemma}
For a map $\iota \colon \aO \to \aO'$ of $\Ninfty$ operads, the
adjunction $\big(\iota_{\redFreeO}^{\redFreeO[\cO']},
\iota_{\redFreeO[\cO']}^{\redFreeO}\big)$ is a Quillen adjunction for the
standard model structures.
\end{lemma}


This allows us to describe a peculiar phenomenon for homotopy
colimits: homotopy colimits in algebras over $\aO$ are always computed
in the underlying category $\GTop_\ast[\redFreeO[O^{tr}]]$, for any
\(\Ninfty\) operad $\aO$.  In particular, for distinct $\Ninfty$
operads $\aO$ and $\aO'$, homotopy colimits in $\GTop_\ast[\redFreeO]$
and $\GTop_\ast[\redFreeO[\bO']]$ coincide.  We can immediately deduce
this from Proposition~\ref{prop:georeal} and the fact that finite
coproducts and products coincide for $\cO$-algebras; we prove this
latter statement as a consequence of an analysis of the free
$\cO$-algebra below in Corollary~\ref{cor:coprodprod}.

\begin{lemma}\label{lem:hocoforg}
Let $\aO$ be an $\Ninfty$ operad.  The forgetful functor $\iota_{\cO}^{\cO^{\tr}}
\colon \GTop_\ast[\redFreeO] \to \GTop_\ast[\redFreeO[\cO^{tr}]]$
creates homotopy colimits. 
\end{lemma}


\subsubsection{Fixed points of \(\cO\)-algebras}
We close this section by recording some standard observations about the
behavior of the categorical fixed-point functor on $\cO$-algebras.

\begin{lemma}\label{lem:FixedPoints}
\hspace{5 pt}
\begin{enumerate}
\item If $B$ is an $\cO$-algebra, then $B^H$ is naturally an $\cO^{tr}$-algebra.

\item If $q\colon G/H\to G/K$ is a $G$-map, then the induced restriction map
\[
q^{\ast}\colon B^{K}\to B^{H}
\]
is naturally a map of $\cO^{tr}$-algebras.
\end{enumerate}
\end{lemma}

\begin{proof}
If $B$ is an $\cO$-algebra, then by restriction, $B$ is an
$\cO^{tr}$-algebra. Since $\cO^{tr}$ has a trivial $G$-action, and
since passage to fixed points is a strong symmetric monoidal on
spaces, $B^{G}$ is automatically an $\cO^{tr}$-algebra.

For the second point, the map $q^{\ast}$ is the $G$-fixed
points of the map 
\[
\Map(G/K,B)\xrightarrow{q^{\ast}}\Map(G/H,B),
\]
which is given by cotensoring out of the map \(G/H\to G/K\), and hence is a map of \(\cO^{tr}\)-algebras.
\end{proof}

This lemma implies that an $\cO$-algebra specifies a diagram over the
orbit category of algebras over $\cO^{tr}$.

\begin{corollary}
There is a contravariant functor 
\[
\GTop_\ast[\redFreeO] \to \Fun\big((\Orb^G)^{\op}, \GTop_\ast[\redFreeO[\cO^{tr}]]\big)
\]
specified by 
\[
B\mapsto (G/H \mapsto B^H) 
\]
and $G/H \to G/K$ goes to $q^{\ast}$.
\end{corollary}

A model for $\cO^{tr}$ is the little disks operad for a
trivial universe $\R^{\infty}$, and so assuming that $B^{H}$ is
grouplike for all $H$, using a functorial non-equivariant delooping
machine we can associate to an object of $\Fun((\Orb^G)^{\op},
\GTop_\ast[\redFreeO[\cO^{tr}]])$ a contravariant functor  
\[
\Orb^G \to\Sp.
\]
This is the data of a naive $G$-spectrum, and passing to homotopy
groups gives rise to a coefficient system of abelian groups.  For
general operads $\cO$, there is more structure on the homotopy groups:

\begin{theorem}[{\cite[7.5]{BHNinfty}}]
If \(B\) is an \(\cO\)-algebra in spaces, then coefficient system
\[
T\mapsto [T_{+}\wedge S^{k},B]=:\m{\pi}_{k}(B)(T)
\]
extends naturally to an \(\cO\)-Mackey functor.
\end{theorem}

Finally, we record two results about the interaction of fixed points
with the tensor and cotensor with (non-equivariant) spaces.

\begin{lemma}\label{lem:tensorfixed}
Let $B$ be an $\cO$-algebra in based $G$-spaces.  For any $H \subset G$ and
based space $X$ that is the geometric realization of a levelwise-finite
simplicial set, the natural map 
\[
X \otimes B^H \to (X \otimes B)^H
\]
is a weak equivalence of $\cO^H$-algebras.
\end{lemma}

\begin{proof}
By Proposition~\ref{prop:georeal}, the tensor with $X$ can be computed
as the geometric realization of the levelwise tensor.  Since $X$ is
levelwise-finite, these tensors are just coproducts.  Therefore, we
can reduce to showing that the natural map 
\[
B^H \vee_{\cO^H} B^H \to (B \vee_{\cO} B)^H
\]
induced by the two inclusions $B \to B \vee B$ is an equivalence for
all $H \subseteq G$.  This now follows from the commutative diagram 
\[
\xymatrix{
B^H \vee_{\cO^H} B^H \ar[d] \ar[r] \ar[d] & (B \vee_{\cO} B)^H \ar[d] \\
B^H \times B^H \ar[r] & (B \times B)^H,
}
\]
as the vertical maps are equivalences by
Corollary~\ref{cor:coprodprod} below, and the bottom horizontal map is an
equivalence because fixed-points and products commute.
\end{proof}

The case of the cotensor is immediate, since cotensor in $\cO^{tr}$ is
formed in the underlying spaces. 

\begin{lemma}
Let $B$ be an $\cO$-algebra in spaces and $X$ a based space.  Then we
have a natural equivalence of $\cO^{G}$-algebras 
\[
(F(X,B))^G\simeq F(X,B^G).
\]
\end{lemma}

\section{The $\cO$-Spanier-Whitehead category and the $\cO$-stable category}

Our model of the equivariant $\cO$-stable category is motivated by
very classical considerations.  The original approach to the stable
category regarded the objects of the Spanier-Whitehead category as
fundamental and the other objects as determined by formal
considerations.

In the Spanier-Whitehead category, the objects are finite based 
CW-complexes, and the morphisms from $X$ to $Y$ are the stable
homotopy classes of maps from $X$ to $Y$, defined as
\[
\{X,Y\} := \colim_n
\big[\Sigma^{n}X,\Sigma^{n}Y\big]\cong
\big[X,\colim_{n}\Omega^{n}\Sigma^{n}Y\big]. 
\]
For connected $X$ and $Y$, we can interpret this description
operadically; May's approximation theorem~\cite[2.7]{gils} implies that 
there is a natural zig-zag of equivalences of $\cO$-algebras
\[
\colim_n \Omega^n \Sigma^n Y \simeq \redFreeO(Y)
\]
for any $E_\infty$ operad $\cO$.  The free-forgetful adjunction then 
provides a natural isomorphism 
\[
\{X,Y\} \cong \big[ \redFreeO(X), \redFreeO(Y)],
\]
where on the right we are working with homotopy classes of maps of
$\cO$-algebras.  One proceeds by showing that the suspension functor
is fully-faithful and then formally inverting it to produce a category
in which suspension is an auto-equivalence.

This definition applies equally well equivariantly.  Here, the objects
are finite based $G$-CW complexes, while the maps are defined as
\[
\{X,Y\}_{U}^{G} := \colim_{V\subset U}
\big[\Sigma^{V}X,\Sigma^{V}Y\big]^{G}\cong \big[X,\colim_{V\subset
    U}\Omega^{V}\Sigma^{V}Y\big]^{G}, 
\]
where the colimit is taken over all finite dimensional subspaces of
some complete universe $U$.  The equivariant analogue of the
approximation theorem (e.g., see~\cite[1.18]{carusowaner}) implies
that for $G$-connected $X$ and $Y$ we have a description
\[
\{X,Y\}_U^{G} \cong [\redFreeO(X), \redFreeO(Y)]_{\GTop_\ast[\redFreeO]}^G,
\]
where $\cO$ is any genuine $G$-$E_\infty$ operad.  Here again, the
suspension by any finite dimensional representation $V$ is
fully-faithful, so we formally invert these.  In fact, by a beautiful
observation of Lewis, at this point it suffices to invert the trivial 
suspensions; all other are automatically inverted in this context.

The operadic description then generalizes immeditely.

\begin{definition}\label{def:OSWCat}
Let $G$ be a finite group and let \(\cO\) be an \(\Ninfty\) operad.  Define the {\defemph{
    $\cO$-Spanier-Whitehead category}} $\SW^{G}_{\cO}$ by
letting the objects be finite based $G$-CW complexes and the
morphisms be 
\begin{align*}
\SW^{G}_{\cO}(X,Y)=\{X,Y\}_{\cO}^{G}&:=\Ho\GTop_\ast[\redFreeO]\big(\redFreeO(\Sigma
X),\redFreeO(\Sigma Y)\big)^{G} \\
&\cong \Ho\GTop_\ast \big(X,\Omega\redFreeO(\Sigma Y)\big)^{G}. 
\end{align*}
Composition is defined via the category $\GTop_\ast[\redFreeO]$.
\end{definition}

As a consequence of Lemma~\ref{lem:groupcompletion} which we prove
below, we obtain the following basic result about the
$\cO$-Spanier-Whitead category.

\begin{theorem}\label{thm:Suspensions}
If $\cO$ is an indexing system for $G$, then the trivial suspension
induces a fully-faithful self-embedding 
\[
\Sigma\colon\SW^{G}_{\cO}\to\SW^{G}_{\cO}.
\]
\end{theorem}

We have an obvious ``stabilization'' functor taking a $G$-space to
an $\cO$-stable object.  

\begin{notation}
Let $\aW_G$ denote the full subcategory of $\GTop$ spanned by the
finite $G$-CW complexes and $\aW_{G,\ast}$ the full subcategory of
$\GTop_\ast$ spanned by the based finite $G$-CW complexes.
\end{notation}

\begin{definition}
Let
\[
\Sigma^{\infty}_{\cO}\colon \aW_{G,\ast} \to \SW^{G}_{\cO}
\]
be the functor specified on objects by the assignment $X \mapsto X$ 
and which takes a morphism \(f\colon X\to Y\) to the induced morphism  
\[
\redFreeO(\Sigma X)\to \redFreeO(\Sigma Y).
\]
Let 
\[
\Sigma^{\infty}_{\cO,+}\colon \aW_G \to \SW^{G}_{\cO}
\]
be the composite of $\Sigma^{\infty}_{\cO}$ with the natural functor
$\Top^{G}\to\Top^{G}_{\ast}$ adding a $G$-fixed disjoint basepoint.
\end{definition}

This description suggests that we regard grouplike \(\cO\)-algebras as
stable objects. Just as classically, these will model connective
spectra. The process of stabilization is then the formal process of
inverting the suspension.

%

\subsection{Equivariant group completions}

A key step in stabilization is group completion.  To study this
process, we work at the level of coefficient systems.
\begin{definition}
\mbox{}
\begin{enumerate}
\item A $\cO$-algebra $B$ is grouplike if $\m\pi_0(B)(T)$ is
a group for every finite \(G\)-set \(T\).

\item A map $A \to B$ of $\cO$-spaces is an equivariant group
  completion if $B$ is grouplike and for every $H \subset G$ the
  induced map of fixed-point spaces $A^H \to B^H$ is a
  (non-equivariant) group completion.
\end{enumerate}
\end{definition}

Given a $\cO$-algebra, we can form the group completion as follows.

\begin{lemma}\label{lem:groupcompletion}
Let $B$ be a cofibrant $\cO$-algebra.  Then the natural map 
\[
B \to F(S^1, S^1 \otimes B)
\]
is an equivariant group-completion.
\end{lemma}

\begin{proof}
For each subgroup $H \subset G$, passing to $H$-fixed points yields
the map
\[
B^H \to F(S^1, S^1 \otimes B)^H.
\]
Since $S^1$ has trivial $H$-action, we can factor this map as the
composite 
\[
B^H \to F(S^1, S^1 \otimes B^H) \to F(S^1, (S^1 \otimes B)^H) \cong
F(S^1, S^1 \otimes B)^H,
\]
and by Lemma~\ref{lem:tensorfixed} all maps except the first one are
weak equivalences.  Moreover, $B^H$ has an action of $\cO^H$, which by
definition is a nonequivariant $E_\infty$ operad.  We know
that if $B$ is cofibrant then $B^H$ has the homotopy type of a
cofibrant algebra over $\cO^H$, by Corollary~\ref{cor:fixcof} below.
Since for any cofibrant algebra $A$ over an $E_\infty$ operad the
natural map 
\[
A \to F(S^1, S^1 \otimes A)
\]
is a group-completion~\cite[6.5]{BasterraMandell}, the result follows.
\end{proof}

To incorporate the group completion into the model category structure,
we proceed as follows.  We can take the map
\[
f \colon \redFreeO S^0 \to F(S^1, S^1 \otimes \redFreeO S^0)
\]
and factor into a cofibration $\tilde{f}$ followed by an acyclic
fibration. Performing Bousfield localization with respect to $\tilde{f}$ then
leads to the following model structure (Compare~\cite[5.3]{SagaveSchlichtkrull}).

\begin{theorem}\label{thm:GroupCompletionModel}
There is a cofibrantly generated $G$-topological model structure, the
{\em group-completion model structure}, on $\GTop_\ast[\redFreeO]$ in which 
\begin{enumerate}
\item the cofibrations are the same as in the standard
model structure
\item the weak equivalences are the maps which induce
weak equivalences after cofibrant-replacement and group-completion, and
\item the fibrations are determined by the right lifting property with
respect to the acyclic cofibrations.
\end{enumerate}
\end{theorem}

The fibrant objects in this category provide our model for connective
\(\cO\)-spectra: In this category, finite wedges are finite products
and the hom sets in the homotopy category are all abelian
groups. Moreover, the homotopy coefficient systems extend naturally to
\(\cO\)-Mackey functors. This is an algebraic shadow of a more general
Wirthm\"uller isomorphism we will prove below.

\subsection{The category of \(\cO\)-spectra}

The category $\GTop_{\ast}[\redFreeO]$, equipped with the group-completion model structure
(or equivalently the full subcategory of group-complete objects in
$\GTop_{\ast}[\redFreeO]$), provides a model of connective $\cO$-spectra. In order to obtain the nonconnective objects, we formally invert the suspension, working with the
category of spectrum objects in $\cO$-spectra. 

\begin{definition}
A spectrum object $E$ in $\GTop_\ast[\redFreeO]$ is a collection of
$\cO$-algebras $\{E_i\}$ for $i \in \bN$ together with natural maps of
\(\cO\)-algebras 
\[
S^1 \otimes E_i \to E_{i+1}
\]
for each $i$.  

A morphism of spectra \(E\to E'\)  is a collection of maps of \(\cO\)-algebras 
\[
E_i \to E'_i
\] 
that commute with the structure maps.  

We write $\Sp^{G}_{\cO}$ to denote the category of spectrum objects and refer to the objects as
$\cO$-spectra.
\end{definition}

\begin{definition}
An $\Omega$-spectrum in $\Sp^{G}_{\cO}$ is defined to be a
spectrum \(E\) such that each adjoint structure map 
\[
E_n \to \Omega E_{n+1}
\]
is a weak equivalence of $\cO$-algebras. 
\end{definition}

In order to construct a model structure on spectrum objects, it is convenient to work with stable
homotopy objects. Since the homotopy coefficient system of an \(\cO\)-algebra is naturally an \(\cO\)-Mackey functor, we have a similar extension.

\begin{definition}
Let $E$ be an $\cO$-spectrum. The homotopy \(\cO\)-Mackey functors of \(E\) are the \(\cO\)-Mackey functors
\[
\m{\pi}_{n}(E):=\colim_{k}\m{\pi}_{n+k}(E_{k}),
\]
where the transition maps in the direct system are given by applying homotopy groups to the adjoint structure maps in the spectrum.
\end{definition}

\begin{definition}
We say a map $E \to E'$ in $\Sp^{G}_{\cO}$ is a stable equivalence if
the induced map of homotopy \(\cO\)-Mackey functors $\m\pi_n(E) \to \m\pi_n(E')$ is an
isomorphism for each $n \in \bZ$.  
\end{definition}

The following lemma is immediate.

\begin{lemma}
A map $f \colon E \to E'$ of $\Omega$-spectra in $\Sp^{G}_{\cO}$ is a
weak equivalence if and only if $E_i \to E'_i$ is a weak equivalence in
$\GTop_{\ast}[\redFreeO]$ for each $i \in \bN$.
\end{lemma}

By the argument for~\cite[3.1]{BasterraMandell}, we can conclude the
following theorem.

\begin{theorem}
There is a model structure on $\Sp^{G}_{\cO}$ where
\begin{enumerate}
\item the cofibrations are the maps \(E\to E'\) such that 
\[
E_0 \to E'_0
\]
is a cofibration and 
\[
(S^1 \otimes E'_n) \smashover{(S^1 \otimes E_n)} E_{n+1} \to E'_{n+1}  
\]
is a cofibration for each $n$,
\item the fibrations are the maps \(E\to E'\) such that
\[
E_n \to E'_n
\]
is a fibration and
\[
E_n \to E'_n \timesover{\Omega E'_{n+1}} \Omega E_{n+1}
\]
is a fibration for each $n$, and 
\item the weak equivalences are the stable equivalences.
\end{enumerate}
(Here we are using the standard model structure on $\GTop[\redFreeO]$, not
the group-complete model structure.)
\end{theorem}

It is straightforward to check directly the following proposition.

\begin{proposition}
A fibrant object in the model structure on $\Sp^{G}_{\cO}$ is an
$\Omega$-spectrum.
\end{proposition}

It suffices to work with the standard model structure when studying
the homotopy theory of $\Sp^{G}_{\cO}$ because for an
$\Omega$-spectrum, $E_0$ is grouplike.  Conversely, given a
$\cO$-algebra $B$, we can form the ``suspension spectrum'' 
\[
(\Sigma^{\infty} B)_n = (S^1 \sma \ldots \sma S^1) \otimes B,  
\]
where the structure maps are given by the canonical maps.  When $B$ is
cofibrant, the adjoint structure maps 
\[
(\Sigma^{\infty} B)_n \to
\Omega (\Sigma^{\infty} B)_{n+1}
\] 
are weak equivalences for $n > 0$
and the group completion for $n = 0$.  In particular, if $B$ is
already grouplike, this is an $\Omega$-spectrum.

\begin{remark}
Given a cofibrant $\cO$-algebra $B$, we can also consider the
stabilization obtained as follows.  Let $\Gamma_G$ denote the category
of finite based $G$-sets, and form the $\Gamma_G$-space specified by
$T \mapsto X \otimes T$.
\end{remark}

We can now describe the suspension spectrum functor from $G$-spaces
\begin{definition}
Let 
\[
\Sigma^{\infty}_{\cO}\colon\Top^{G}_{\ast}\to \Sp^{G}_{\cO}
\]
be the functor
\[
X\mapsto \Sigma^{\infty} \big(\redFreeO(X)\big).
\]

Let
\[
\Sigma^{\infty}_{\cO,+}\colon\Top^{G}\to \Sp^{G}_{\cO}
\]
be the composite of \(\Sigma^{\infty}_{\cO}\) with the functor $(-)_+$
adding a disjoint basepoint.
\end{definition}

We note that the suspension spectrum of a $G$-space is always an
$\Omega$-prespectrum above level $0$, by
Lemma~\ref{lem:groupcompletion}.

\begin{corollary}\label{cor:omega}
Let $X$ be a based and connected $G$-CW complex.  Then the suspension
spectrum $\Sigma^{\infty}_{\cO} X$ is an $\Omega$-spectrum except at
level $0$, where it is the group completion.
\end{corollary}

\subsection{Structural properties of $\cO$-spectra}

We now turn to the structural properties of the categories
$\Sp^{G}_{\cO}$.  

\subsubsection{Induction, coindunction, and Wirthm\"uller isomorphisms}

We begin by enumerating the construction of the usual change of group
and fixed-point functors.  For $H \subset G$, there are natural
forgetful functors.

\begin{definition}
For each $H \subset G$, there is a forgetful functor 
\[
\iota^*_H \colon \Sp^{G}_{\cO} \to \Sp^{H}_{i_H^\ast\cO}
\]
specified by the levelwise application of the forgetful functor on
$\cO$-algebras 
\[
\iota^*_H \colon \GTop_{\ast}[\redFreeO] \to \Top^{H}_{\ast}[\iota^*_H \redFreeO].
\]
\end{definition}

The forgetful functors have homotopical left and right adjoints. The
right adjoint is easier to describe, since it takes place entirely in
the underlying category of \(G\)-spaces.

\begin{proposition}[See e.g., {\cite[\S 6.2.2]{BHNinfty}}]
If \(B\) is an \(i_{H}^{\ast}\cO\)-algebra in based \(H\)-spaces, then
\(\Map^{H}(G_{+},B)\) is naturally a \(\cO\)-algebra in \(G\)-spaces
and the coinduction functor 
\[
B\mapsto \Map^{H}(G_{+},B)
\]
is right adjoint to the forgetful functor
\(i_{H}^{\ast}\).  The adjunction is a Quillen adjunction.
\end{proposition}

\begin{proof}
Let $B$ be an $\i_H^\ast \cO$-algebra.  We produce the action map
\[
\cO(n) \times (\Map^H(G_+, B))^{\times n} \to \Map^H(G_+, B)
\]
as the adjoint of the composite 
\begin{align*}
i_H^\ast (\cO(n) \times (\Map^H(G_+, B))^{\times n}) &\cong i_H^\ast
\cO(n) \times (i_H^\ast \Map^H(G_+, B))^{\times n} \\ &\to i_H^\ast \cO(n)
\times B^{\times n} \to B,
\end{align*}
where we are using the fact that $\Map^H(G_+,-)$ is the adjoint to
$i_H^\ast$ in spaces.  It is straightforward to check that these
action maps satisfy the required compatibilities, and the fact that
this is the right adjoint in operadic algebras is a consequence of the
construction.  To see that this is a Quillen adjunction, observe that
$\Map^H(G_+,-)$ preserves fibrations since these are determined by the
forgetful functor to $G$-spaces.
\end{proof}

We can extend the right adjoint to spectra by applying it levelwise,
getting a functor 
\[
F^{H}(G_{+},-) \colon \Sp^H_{i^H_\ast \cO} \to \Sp^G_{\cO}.
\]

\begin{corollary}
The functor
\[
E\mapsto F^{H}(G_{+},E)
\]
is right adjoint to \(i_{H}^{\ast}\) and this is a Quillen adjunction.
If \(E\) is an \(\Omega\)-\(i_{H}^{\ast}\cO\)-spectrum, then
\(F^{H}(G_{+},E)\) is naturally an \(\Omega\)-\(\cO\)-spectrum.
\end{corollary}

The left adjoint is slightly harder to describe, since it is built out
of the operadic coproduct, rather than the underlying one. Again, we
begin in \(\cO\)-algebras in spaces. 

\begin{proposition}
The forgetful functor \(i_{H}^{\ast}\) has a left adjoint 
\[
G_{+}\otimes_{H}(\mhyphen)\colon\Top^{H}[i_{H}^{\ast}\bO]\to\Top^{G}[\bO].
\]
\end{proposition} 

The construction is fairly straightforward. By the universal property,
on free objects, we have 
\[
G_{+}\otimes_{H}\tilde{\mathbb P}_{i_{H}^{\ast}\cO}(X)\cong \redFreeO(G_{+}\wedge_{H}X).
\]
Writing any general operadic algebra as a reflexive coequalizers (or
more generally a sifted colimit like a geometric realization) yields the standard formula for the left adjoint in general.

Applying this levelwise gives us a functor 
\[
G_{+}\otimes_{H}(\mhyphen)\colon\Sp^{H}_{i_{H}^{\ast}\cO}\to\Sp^{G}_{\cO}.
\]

\begin{corollary}
The functor \(G_{+}\otimes_{H}(\mhyphen)\) is left adjoint to
\(i_{H}^{\ast}\).  The adjunction is a Quillen adjunction.
\end{corollary}

\begin{proof}
To see that this is a Quillen adjunction, observe that it is immediate
that $i_H^\ast$ preserves weak equivalences and fibrations.
\end{proof}

A defining property of the genuine equivariant category is that
induction and coinduction agree.  In the setting of $\cO$-spectra,
this is true for admissible sets.  There is a canonical natural
transformation $H_{+}\otimes_{K}(\mhyphen)\Rightarrow
F^{K}(H_{+},\mhyphen)$, and we prove the following comparison result
below as Theorem~\ref{thm:OWirth}.

\begin{theorem}
If \(H/K\) is an admissible \(H\)-set for \(\cO\), then the natural
transformation 
\[
H_{+}\otimes_{K}(\mhyphen)\Rightarrow F^{K}(H_{+},\mhyphen)
\]
induces a natural isomorphism in \(\Ho(\Sp^{H}_{i_{H}^{\ast}\cO})\).
\end{theorem}

These equivalences endow the set of homotopy classes of maps between
$\cO$-spectra with extra structure; they are automatically coefficient
systems, and the Wirthm\"uller isomorphisms endows them with
transfers.

\begin{corollary}
If \(E\) and \(E'\) are spectra, then the coefficient system
\[
T\mapsto [T_{+}\wedge E,E']
\]
extends naturally to an \(\cO\)-Mackey functor.
\end{corollary}
All algebraic invariants represented by \(\cO\)-spectra naturally then
give \(\cO\)-Mackey functors.

The question of which Wirthm\"uller isomorphisms we have is closely
connected to various kinds of rigidity of compact generators in the
homotopy category of \(\cO\)-spectra. Since we can build cofibrant
\(\cO\)-spectra out of free \(\cO\)-algebras on \(G\)-CW complexes, we
immediately deduce that the compact generators are again suspension
spectra of orbits. 

\begin{proposition}\label{prop:generators}
The category \(\Ho(\Sp^{G}_{\cO})\) has a set of compact generators
given by the suspension spectra of orbits
\(\Sigma^{\infty}_{\cO,+}G/H\).
\end{proposition}

\begin{remark}
Although we have not constructed the symmetric monoidal structure on
the category $\Ho(\Sp^{G}_{\cO})$, in principle the Wirthm\"uller
isomorphism then tells us when the compact generators are dualizable.
Specifically, we deduce from the Wirthm\"uller isomorphism applied to
\(E=S^{0}\) that if \(H/K\) is an admissible \(H\)-set, then \(H/K\)
is dualizable (and its dual is itself).  That other orbits are not
self-dual is an easy consequence of the tom Dieck splitting
established below.  In particular, if we are not in the genuine case,
then the associated tensor triangulated category is not rigid in the
sense of \cite{BalmerPic}.
\end{remark}

\subsection{Fixed and geometric fixed points}
The usual fixed-point constructions are quite sensitive to the operad that shows up. We begin by discussing change of operad in the context of the
$\cO$-stable category; the results here are a straightforward
extension of the results in the context of operadic algebras in
spaces.

\begin{proposition}
Let $\cO$ and $\cO'$ be $\Ninfty$ operads such that there is a map 
$\cO \to \cO'$.  Then there are adjoint functors 
\[
\iota_{\cO}^{\cO'} \colon \Sp^{G}_{\cO} \to \Sp^{G}_{\cO'}
\]
and
\[
\iota_{\cO'}^{\cO} \colon \Sp^{G}_{\cO'} \to \Sp^{G}_{\cO}
\]
that form a Quillen adjunction.
\end{proposition}

The $\cO$-stable category is entirely determined by the indexing
system of $\cO$ in the following sense.  

\begin{corollary}
If $\cO \to \cO'$ is a weak equivalence of $\Ninfty$ operads, then
the adjunction $(\iota_{\cO}^{\cO'}, \iota_{\cO'}^{\cO})$ is a Quillen
equivalence.
\end{corollary}

\begin{proof}
This follows from the fact that we can compute $\iota_{\cO'}^{\cO} X$
using the bar construction $B(\cO', \cO, X)$; if $\cO \to \cO'$ is a
weak equivalence, this bar construction is equivalent to $B(\cO, \cO,
X)$ which is always homotopic to $X$.
\end{proof}

\subsubsection{Categorical fixed points}
There are categorical fixed points functors which are adjoint to the appropriate ``trivial'' spectra. Let \(N\) be a normal subgroup of \(G\). Then \(\cO^{N}\) is an \(\Ninfty\)-operad for both \(G\) and for \(Q=G/N\), and the obvious inclusion 
\[
\cO^{N}\hookrightarrow\cO
\]
is a map of \(G\) \(\Ninfty\)-operads. Since \(N\) acts trivially on \(\cO^{N}\), the inclusion of \(G/N\)-spaces into \(G\)-spaces gives an inclusion
\[
q^{\ast}\colon \Sp^{G/N}_{\cO^{N}}\to\Sp^{G}_{\cO^{N}}.
\]
\begin{definition}
The pushforward from \(G/N\) \(\cO^{N}\)-spectra to \(G\) \(\cO\)-spectra is the composite
\[
\iota_{\mathbb P_{\cO^{N}}}^{\mathbb P_{\cO}}\circ q^{\ast}\colon\Sp^{G/N}_{\cO^{N}}\to\Sp^{G}_{\cO}.
\]
\end{definition}

Just as classically, the pushforward has a right adjoint.

\begin{definition}
The categorical \(N\)-fixed-point functor
\[
(-)^N \colon \Sp^{G}_{\cO} \to \Sp^{G/N}_{\cO^{N}}
\]
is specified by the space-level $N$-fixed points functor applied to the restriction to an \(\cO^{N}\)-spectrum.
\end{definition}

Unraveling the adjunctions gives the following, just as classically.

\begin{proposition}
The \(N\)-fixed points functor is the right adjoint to $\iota_{\mathbb
  P_{\cO^{N}}}^{\mathbb P_{\cO}}\circ q^{\ast}$, and this is a Quillen
adjunction.
\end{proposition}

There is a version of the tom Dieck splitting, allowing us to identify the fixed points of suspension spectra. Since \(\cO^{G}\) is just an ordinary \(E_{\infty}\) operad, a group-like \(\cO^{G}\)-algebra is a connective spectrum. Our analysis of the fixed points of free algebras in Section~\ref{sec:tomd} below gives a proof of this \(\cO\)-algebra tom Dieck splitting.

\begin{theorem}[$\cO$-Spectra tom Dieck]\label{thm:OsptomDieck}
Let $X$ be a $G$-CW complex. Then we have a natural equivalence 
\[
\big(\Sigma^{\infty}_{\cO} X_+\big)^G\simeq
\bigvee_{G/H\in\pi_0\cO(G)} 
EW_G(H)_{+}\smashover{W_G(H)} \Sigma^{\infty}_+X^H, 
\]
where here the index for the wedge is all isomorphism classes of
orbits \(G/H\) which are admissible for \(\cO\). 
\end{theorem}

\subsubsection{Geometric fixed points}

For many purposes it is desirable to have a fixed-point functor which
has the property that the value on $\Sigma^{\infty}_{\cO} X_+$ would
be $\Sigma^{\infty} X^G$.  This is accomplished by the geometric fixed
points.  Heuristically, the geometric fixed points remove the ``algebraic''
part of the fixed points arising from transfer maps, leaving only the
last piece of the terms in the tom Dieck splitting.

Homotopically, this is achieved by the composite of two functors:
\begin{enumerate}
\item the Bousfield localization which nullifies the localizing,
  triangulated subcategory generated by the induced cells \(G/H_{+}\)
  for \(H\) a proper subgroup, followed by 
\item the categorical fixed points, which induces an equivalence of
  categories between the local objects and the stable category.
\end{enumerate}

We can carry this for $\cO$-spectra in exactly the same way it is done
classically.

\begin{definition}
The homotopical geometric fixed-point functor
\[
\Phi^G \colon \Sp^{G}_{\cO} \to \Sp
\]
is specified by the formula 
\[
E \mapsto ((\tilde{E}\cP \sma E)^{\fib})^G.
\]
\end{definition}

By construction, $\Phi^G (-)$ preserves weak equivalences between
cofibrant objects.  The ``localization'' part is given by tensoring
with \(\tilde{E}\cP\).  By definition, the restriction of this to any
proper subgroup is equivariantly contractible, so the result is also
true for the tensored \(\tilde{E}\cP\wedge E\). We then compose with
the categorical fixed points.  There is always a canonical map
\[
E^G \to \Phi^G E.
\]

Note that in contrast to the situation with traditional models of
spectra, we do not have a good point-set construction of the geometric
fixed points.

\begin{remark}
For a general \(\cO\), the canonical map
\[
E^{G}\to \Phi^{G}(E)
\]
can be closer to an equivalence than expected from the usual genuine
case.  Notably, if the only admissible \(G\)-sets are those with a
trivial action, then this is an equivalence.  This holds, for example,
for \(\cO^{tr}\) (which is to say that the fixed points of naive
\(G\)-spectra agree with the geometric fixed points), but it also
holds for the \(C_{4}\)-linear isometries operads associated to
\(\infty(\mathbb R\oplus \lambda)\) described above. 
\end{remark}

We turn now to proving the key structural properties of
\(\cO\)-spectra, showing the tom Dieck splitting and the Wirthm\"uller
isomorphism.  

\section{The tom Dieck splitting of a free $\cO$-algebras}\label{sec:tomd}

Since our model of the $\cO$-stable category is built from the
category of $\cO$-algebras in $\GTop$, the core of our analysis of the
category of $\cO$-spectra will depend on describing the equivariant
homotopy type of free $\cO$-algebras.  The purpose of this section is
to calculate this homotopy type.  In particular, we establish a
splitting theorem that is a version of the tom Dieck splitting in the
category of $\cO$-algebras.

\subsection{Fixed points of free $\cO$-algebras}

The proof of the following theorem is the goal of this subsection.  In
the statement, we write $\Map^G(-,-)$ to denote the set of $G$-maps
between $G$-spaces and $\Aut_G(T)$ for the automorphisms of \(T\) as a \(G\)-set.  Note that since
$\Aut_G(T)$ is a group under composition (with unit the identity
map), we can form the universal space $E \Aut_G(T)$ via the two-sided
bar construction  
\[
E \Aut_G(T) = B(\Aut_G(T), \Aut_G(T), *).
\]
For an $\Ninfty$ operad $\aO$, we write $\pi_0\cO(G)$ to denote the
set of isomorphism classes of $G$-sets which are admissible for $\cO$. 

\begin{theorem}\label{thm:OtomDieck}
Let $\cO$ be an $\Ninfty$-operad, and let $X$ be a $G$-CW
complex. Then for all $n$ we have a natural weak equivalence 
\[
\big(\cO_n\timesover{\Sigma_n} X^n\big)^G\simeq
\coprod_{\substack{T\in\pi_0\cO(G),\\ |T|=n}}
E\Aut_G(T)\timesover{\Aut_G(T)} \Map^G(T,X). 
\]
Moreover, if $\cO\to\cO'$ is a map of $\Ninfty$ operads, then the
natural map  
\[
\big(\cO_n\timesover{\Sigma_n}X^n\big)^G\to \big(\cO'_n\timesover{\Sigma_n} X^n\big)^G
\]
is homotopic to the obvious inclusion of wedge summands induced by the
inclusion $\pi_0 \cO(G)\subset \pi_0 \cO'(G)$. 
\end{theorem}

As an immediate consequence, we get the following description of the
fixed-points of the free $\aO$-algebra on a space $X$.

\begin{corollary}\label{cor:FreeOAlg}
Let $\cO$ be an $\Ninfty$-operad, and let $X$ be a $G$-CW
complex. Then we have a natural weak equivalence 
\[
\big(\FreeO(X)\big)^G \simeq\coprod_{T\in\pi_0\cO(G)}
E\Aut_G(T)\timesover{\Aut_G(T)} \Map^G(T,X). 
\]
\end{corollary}

Throughout the course of this section, we will refine this corollary
several times, deducing increasingly strong statements.

\begin{remark}
When $\cO$ is the trivial $\Ninfty$ operad, then there is a unique
isomorphism class of admissible $G$-sets of cardinality $n$, namely
the set of $n$ elements with the trivial action.  In this case, we
recover the classical results
\[
\big(\cO^{tr}_n\timesover{\Sigma_n} X^n\big)^G\simeq E\Sigma_n\timesover{\Sigma_n} (X^G)^n,
\]
and
\[
\big(\mathbb P_{\cO^{tr}}(X)\big)^G\simeq \coprod_{n\geq 0}E\Sigma_n\timesover{\Sigma_n}
(X^G)^n\simeq\mathbb P_{i_e^{\ast}\cO^{tr}}(X^G). 
\]
In words, the fixed points of the free $\aO$-algebra over a trivial
$\Ninfty$ operad on $X$ is the free $E_\infty$ algebra on $X^G$.
\end{remark}

Before we begin, we need a small lemma which helps in working with
universal spaces and determining fixed points thereof.  We first
record a categorical observation about families of subgroups.  Recall
that a sieve in a category $\aC$ is a subcategory $\aS$ such that for
any object $X$ in $\aC$, all morphisms $Z \to X$ are in $\aS$.

\begin{lemma}[{(e.g., \cite[3.5]{BHNinfty},
      \cite[6.3]{GuillouMayRubin})}]\label{lem:sieve}
Let $G$ be a finite group.  Families of subgroups of $G$ are
equivalent to sieves in the orbit category $\cOrb^G$. 
\end{lemma}

Because of this, we will tacitly elide the distinction in the
following discussion, identifying a family $\cF$ with the sieve in
$\cOrb^G$.

We now develop some machinery in order to use the model of the
classifying space of a family provided by Elmendorf's theorem.  Recall
that an $\cOrb^G$-space $X$ is a functor $(\cOrb^G)^{op} \to \Top$.
Given a covariant functor $\cG \colon \cOrb^G \to \Top^G$, we can form
the the two-sided bar construction $B_\cdot(X, \cOrb^G, \cG)$,
which is the simplicial object in $G$-spaces with $k$-simplices
specified by the assignment 
\[
[k] \mapsto \coprod_{G/H_k\to\dots\to G/H_0} X(G/H_0) \times \cG(G/H_k),
\]
where $G$ acts via the action on $\cG(G/H_k)$.  The geometric
realization then becomes a $G$-space, as the action is evidently
compatible with the simplicial structure maps.

We now recall two particularly important examples of inputs for this
two-sided bar construction.

\begin{definition}[{\cite{ElmendorfTheorem}}]
Let $\cF$ be family of subgroups of $G$. Let
$\cF^{op}\colon(\cOrb^G)^{op}\to \Top$ be the functor defined by 
\[
\cF^{op}(G/H)=\begin{cases}
\ast & H\in\cF \\
\emptyset & H\notin\cF.
\end{cases}
\]
Let $J\colon\cOrb^G\to\Top^G$ be the functor which sends an orbit to
itself (viewed now as a $G$-space). 
\end{definition}

The following lemma gives a convenient simplification of this
two-sided bar construction when $X$ is the functor $\cF^{op}$, using
the characterization of Lemma~\ref{lem:sieve}.

\begin{lemma}\label{lem:UniversalSpaces}
Let $\cF$ be a family of subgroups of a group $G$, and let $\cC$ be
any covariant functor $\cOrb^G\to\Top^G$. Then we have an equivariant
isomorphism of simplicial $G$-spaces
\[
B_\bullet(\cF^{op}, \cOrb^G, \cC)\cong B_\bullet(\ast, \cF, \cC|_{\cF}\big),
\]
where $\ast$ is the constant functor from $\cF$ with value a point.
\end{lemma}

\begin{proof}
Recall that a point in $B_k(\cF^{op}, \cOrb^G, \cC)$ is a $(k+2)$-tuple
\[
\big(x\in \cF^{op}(G/H_0), G/H_k\to\dots\to G/H_0, c\in\cC(G/H_k)\big),
\]
and $G$ acts on this via the action in the last coordinate.  Since
$\cF^{op}$ is non-empty only on the sieve $\cF$, any such point
necessarily has $G/H_0\in\cF$, and since $\cF$ is a sieve, this
implies that all of $G/H_0,\dots, G/H_k$ are in $\cF$.  That is, $x$
is in fact a point in $B_k(\ast, \cF, \cC|_{\cF})$.  This
identification is clearly compatible with the simplicial structure
maps.
\end{proof}

We are now ready to perform the basic calculation of the fixed-points
of the extended powers for an $\Ninfty$ operad.  In the following
statement, recall that an admissible set $T$ in an indexing system
$\aO$ corresponds to a homomorphism $\Gamma \colon H \to \Sigma_{|T|}$
such that the associated family contains the graph of $\Gamma$.  We
will refer to the subgroup associated to $T$ as $\Gamma_T$.

\begin{theorem}\label{thm:OHomotopyOrbits}
Let $\cO$ be an $\Ninfty$ operad. Then for all $n$ and for all
$G\times\Sigma_n$-CW complexes $Y$, we have a natural equivalence 
\begin{multline*}
(\cO_n\timesover{\Sigma_n} Y)^G\simeq \\
\coprod_{\substack{T\in\pi_0\cO(G), \\ |T|=n}} E\Aut_{G\times\Sigma_n}(G\times\Sigma_n/\Gamma_T)\timesover{\Aut_{G\times\Sigma_n}(G\times\Sigma_n/\Gamma_T)} \big((G\times\Sigma_n/\Gamma_T)\timesover{\Sigma_n} Y\big)^G.
\end{multline*}

Moreover, if $\cO\to\cO'$ is a map of $\Ninfty$ operads, then the natural map
\[
(\cO_n\timesover{\Sigma_n}Y)^G\to(\cO'_n\timesover{\Sigma_n} Y)^G
\]
is homotopic to the inclusion of summands induced by $\cO(G)\subset\cO'(G)$.
\end{theorem}

\begin{proof}
By hypothesis, $\cO_n$ is a universal space for the family $\cF_n$ of
subgroups of $G\times\Sigma_n$ of the form $\Gamma_T$, where $T$
ranges over all admissible sets of cardinality $n$ for all subgroups.
Moreover, $\cO_n$ is $\Sigma_n$-free and has the homotopy type of
a CW-complex; we can replace $\cO_n$ up to weak equivalence, and
therefore we have flexibility to use any point-set model we find
convenient for this homotopy type.  For our purposes, we will use
Elmendorf's bar construction~\cite{ElmendorfTheorem} as simplified by
Lemma~\ref{lem:UniversalSpaces}.

In particular, we have for any $G\times\Sigma_n$-CW complex $Y$ a
canonical equivariant homeomorphism
\begin{align*}
E\cF_n\times Y\cong 
\big|B_\bullet(\cF_n^{op}, \cOrb^G, J)\big| \times Y &\cong 
\big|B_\bullet(\ast, \cF_n, J|_{\cF_n})\big| \times Y
\\
&\cong \big|B_\bullet(\cF_n^{op}, \cF_n, J|_{\cF_n}\times
Y)\big|,
\end{align*}
where $J|_{\cF_n}\times Y$ is the functor which sends an orbit
$G\times\Sigma_n/\Lambda$ to its product with $Y$; here we are using
the fact that finite products commute with geometric realization.
Furthermore, we have a $G$-equivariant homeomorphism
\[
E\cF_n\timesover{\Sigma_n} Y\cong \big|B_\bullet(\ast, \cF_n, J|_{\cF_n}\timesover{\Sigma_n} Y)\big|,
\]
since colimits commute.  Using the fact that the geometric realization
can be described as a sequential colimit of pushouts along closed
inclusions, we can conclude that fixed points commute with geometric
realization: 
\[
(E\cF_n\timesover{\Sigma_n} Y)^G\cong \Big|B_\bullet\big(\ast, 
\cF_n, (J|_{\cF_n}\timesover{\Sigma_n} Y)^G\big)\Big|. 
\]

If $\Lambda\in \cF_n$, then $\Lambda$ is of the form $\Gamma_T$ for some subgroup $H\subset G$ and some admissible finite $H$-set of cardinality $n$. If $H\subsetneq G$, then in fact $\Gamma_T\subset H\times\Sigma_n$, and we have an isomorphism of $G\times\Sigma_n$-spaces
\[
G\times\Sigma_n/\Gamma_T\cong G\times_H\big(H\times\Sigma_n/\Gamma_T\big).
\]
We therefore deduce that for any $G\times\Sigma_n$-space $Y$, we have a $G$-equivariant homeomorphism
\[
\big(G\times\Sigma_n/\Gamma_T\big)\timesover{\Sigma_n} Y\cong G\times_H \Big(\big(H\times\Sigma_n/\Gamma_T\big)\timesover{\Sigma_n} i_{H\times\Sigma_n}^{\ast} Y\Big).
\]
In particular, if $H\subsetneq G$, then this has no $G$-fixed points. 

Putting this together, we see that the fixed points
$(E\cF_n\timesover{\Sigma_n} Y)^G$ can be computed as the geometric
realization of a simplicial object with $k$-simplices: 
\begin{multline*}
\left\{
\Big((G\times\Sigma_n/\Gamma_{T_k})\to\dots\to (G\times\Sigma_n/\Gamma_{T_0}), y\in \big((G\times\Sigma_n/\Gamma_{T_k})\timesover{\Sigma_n} Y \big)^G\Big)  \right.\\ \Big. \mid
T_0,\dots T_k \text{ admissible for }\cO, T_k\in\cO(G)
\Big\}
\end{multline*}
Now if $T$ is a $G$-set, then $\Gamma_T$ is a maximal element in the
family of subgroups. In particular, all of the maps in the
$k$-simplices above must be isomorphisms.  

Thus, when we pass to fixed points, we are simply restricting to the
maximal subgroupoid of $\cF_n$ containing the maximal subgroups: 
\begin{multline*}
(E\cF_n\timesover{\Sigma_n} Y)^G\cong \\
\coprod_{\substack{T\in\pi_0\cO(G),\\ |T|=n}} E\Aut_{G\times\Sigma_n}(G\times\Sigma_n/\Gamma_T)\timesover{\Aut_{G\times\Sigma_n}(G\times\Sigma_n/\Gamma_T)} \big((G\times\Sigma_n/\Gamma_T)\timesover{\Sigma_n} Y\big)^G.
\end{multline*}

For the second part of the statement, note that the map of $\Ninfty$
operads $\cO\to\cO'$ induces a map on $n$th spaces. Since these are
universal spaces, there is a unique homotopy class of such a map, so
we can without loss of generality assume that the map is the one
induced by the inclusion of the family of $\cO_n$ into that for
$\cO_n'$.
\end{proof}

We can now restrict to the case that $Y=X^n$ for some $G$-space $X$
with the evident induced $G\times\Sigma_n$-action.  Here the orbits in
Theorem~\ref{thm:OHomotopyOrbits} can be put in a more illuminating form.

\begin{lemma}[{\cite[6.2]{BHNinfty}}]\label{lem:OrbitsandCoInduction}
For all $G$-spaces $X$ and all finite $G$-sets $T$ of cardinality $n$,
we have a natural $G$-equivariant homeomorphism
\[
(G\times\Sigma_n/\Gamma_T)\timesover{\Sigma_n} X^n\cong \Map(T,X),
\]
where the mapping space has the conjugation action.
\end{lemma}

Next, we can identify the Weyl group of $\Gamma_T$ as a subgroup of $G \times
\Sigma_n$.

\begin{lemma}\label{lem:OWeylAction}
Let $T$ be a finite $G$-set of cardinality $n$, and let $X$ be a $G$-space. 
\begin{enumerate}
\item We have a natural isomorphism
\[
\Aut_{G\times\Sigma_n}(G\times\Sigma_n/\Gamma_T)=W_{G\times\Sigma_n}(\Gamma_T)\cong \Aut_G(T).
\]
\item The isomorphism of Lemma~\ref{lem:OrbitsandCoInduction} is
  $\Aut_G(T)$ equivariant, where the lefthand side has the evident
  action induced by the isomorphism above and $\Map(T,X)$ is endowed
  with the obvious action of $\Aut_G(T)$.
\end{enumerate}
\end{lemma}

\begin{proof}
We begin by identifying the normalizer of $\Gamma_T$.  Let $f$ be the
homomorphism $G\to\Sigma_n$ defining $T$; the graph of $f$ is
precisely $\Gamma_T$.  Then the normalizer of $\Gamma_T$ in
$G\times\Sigma_n$ is 
\[
N_{G\times\Sigma_n}(\Gamma_T)=\{
(g,\sigma) \mid \forall h\in G, (g,\sigma)\big(h, f(h)\big)(g^{-1}, \sigma^{-1})\in\Gamma_T
\}.
\]
Unpacking the condition, we see that $(g,\sigma)\in N_{G\times\Sigma_n}(\Gamma_T)$ if and only if for all $h\in G$, 
\begin{equation}\label{eq:form}
\sigma f(h)\sigma^{-1}=f(g) f(h) f(g)^{-1}.
\end{equation}
Therefore, we can define a natural homomorphism 
\[
\phi\colon\Aut_G(T)\to N_{G\times\Sigma_n}(\Gamma_T)\to
W_{G\times\Sigma_n}(\Gamma_T) 
\]
as the composite of the map which sends an automorphism $\sigma$ to
the pair $(e,\sigma)$ followed by the projection.

We now deduce two consequences of equation~\eqref{eq:form}:
\begin{itemize}
\item Taking $h=g$, we have
\begin{equation}\label{eq1}
\sigma f(g)\sigma^{-1}=f(g) f(g) f(g)^{-1} = f(g),
\end{equation}
i.e., $f(g)$ and $\sigma$ commute.
\item Multiplying by $\sigma$ on the right and $f(g)^{-1}$ on the
  left, we obtain
\begin{equation}\label{eq2}
f(g)^{-1} \sigma f(h) = f(h) f(g)^{-1} \sigma, 
\end{equation}
i.e., $f(g)^{-1}\sigma$ commutes with $f$ and thus specifies an
automorphism of $T$.
\end{itemize}

We now use these facts to show that $\phi$ is an isomorphism.  Recall
that for any $g\in G$, by definition the element
$\big(g,f(g)\big)\in\Gamma_T$.  The following equalities
\[
(g,\sigma)\simeq (g,\sigma)\big(g^{-1},f(g)^{-1}\big)=(e,\sigma f(g)^{-1})=(e,f(g)^{-1}\sigma)=\phi\big(f(g)^{-1}\sigma\big),
\]
now show that any element $(g,\sigma) \in N_{G \times
  \Sigma_n}(\Gamma_T)$ is equivalent modulo $\Gamma_T$ to something in
the image of $\phi$.  That is, $\phi$ is surjective.  Here the second
equality uses equation~\eqref{eq1} and the third equality uses
equation~\eqref{eq2}.  Similarly, if $\phi(\sigma)=\phi(\sigma')$,
then there is some $g\in G$ such that
\[
(e,\sigma)\big(g, f(g)\big)=(e,\sigma').
\]
This forces $g=e=f(g)$, and hence $\sigma=\sigma'$; $\phi$ is
injective.

The second part is now obvious: passage to $\Sigma_n$ orbits allows us
to move a $G$-equivariant automorphism $\sigma$ of $T$ (viewed as a
particular element in $\Sigma_n$) to the exponent of $X^n$, which is
precisely the definition of the action of $\Aut_G(T)$ on the set of $n$
elements.
\end{proof}

We now put this all together.

\begin{proof}[Proof of Theorem~\ref{thm:OtomDieck}]
Applying theorem~\ref{thm:OHomotopyOrbits} to $Y=X^n$ yields the
formula 
\begin{multline*}
(\cO_n\timesover{\Sigma_n} X^n)^G\cong \\
\coprod_{\substack{T\in\pi_0\cO(G),\\ |T|=n}} E\Aut_{G\times\Sigma_n}(G\times\Sigma_n/\Gamma_T)\timesover{\Aut_{G\times\Sigma_n}(G\times\Sigma_n/\Gamma_T)} \big((G\times\Sigma_n/\Gamma_T)\timesover{\Sigma_n} X^n\big)^G.
\end{multline*}
Using lemma~\ref{lem:OrbitsandCoInduction}, we identify the fixed
points of the orbits in the summand associated to $T$:  
\[
\big((G\times\Sigma_n/\Gamma_T)\timesover{\Sigma_n} X^n\big)^G\cong\Map^G(T,X).
\]
Finally, Lemma~\ref{lem:OWeylAction} identifies the group
$\Aut_{G\times\Sigma_n}(G\times\Sigma_n/\Gamma_T)$ and the action.
\end{proof}

\subsection{The tom Dieck splitting for $\cO$-algebras}

We now apply the identifications of the previous section to deduce a
tom Dieck splitting for $\cO$-algebras.  Specifcally, we will apply
the identification of Corollary~\ref{cor:FreeOAlg} and the naturality
in the operad $\cO$.

\subsubsection{The unbased case}

We begin with some notation and an elementary observation.

\begin{notation}
Let $(H)_{G}$ denote the set of $G$-conjugacy classes of $H$. If $G$ is clear from the context, we will also denote this simply as $(H)$.
\end{notation}

Since an automorphism preserves stabilizers of points, no automorphism
can change orbit types:

\begin{proposition}
If 
\[
T=\coprod_{(H)} n_H\cdot G/H,
\]
where the coproduct ranges over conjugacy classes of subgroups, then we have an identification
\[
\Aut_G(T)\cong \prod_{(H)} \Aut_G(n_H\cdot G/H).
\]
\end{proposition}

Coupled with the universal property of the product, we have the
following further refinement of Corollary~\ref{cor:FreeOAlg}.  

\begin{proposition}\label{prop:FixedofFrees}
Let $\cO$ be an $\Ninfty$ operad, and let $X$ be a $G$-CW
complex. Then we have an identification 
\[
\big(\FreeO(X)\big)^G\simeq \prod_{G/H\in \pi_0 \cO(G)} \coprod_{n\geq 0} E\Aut_G(n\cdot G/H)\timesover{\Aut_G(n\cdot G/H)} (X^H)^n.
\]

If $\cO\subset\cO'$, then the identification can be chosen as the
inclusion of those factors corresponding to $G/H$ with $G/H\in
\pi_0\cO(G)$. 
\end{proposition}

\begin{proof}
Since any finite $G$-set is a disjoint union of orbits
and a finite set is admissible if and only if all of the orbits
occurring in such a decomposition are admissible, the identification
is essentially a rewriting of Corollary~\ref{cor:FreeOAlg}.  The only
further refinement is the application of the isomorphism
\[
\Map^G(n\cdot G/H,X)\cong \Map^G(G/H,X)^n\cong (X^H)^n.
\]
\end{proof}

Next, we record the standard computation of the $G$-automorphisms of
the $G$-set 
\[
n\cdot G/H = \coprod_n G/H.
\]

\begin{proposition}\label{prop:AutnGH}
Let $H\subset G$.  Then we have an isomorphism
\[
\Aut_G(n\cdot G/H)\cong \Aut_G(G/H)\wr \Sigma_n\cong
W_G(H)\wr\Sigma_n.
\]
\end{proposition}

We are now ready to prove the main theorem of this subsection.

\begin{theorem}\label{thm:BPQ}
Let $\cO$ be an $\Ninfty$ operad, and let $X$ be a $G$-CW complex. Let
$\mathcal E$ denote any non-equivariant $E_\infty$ operad. Then we
have a decomposition of $\mathcal E$-algebras 
\[
\big(\FreeO(X)\big)^G\simeq \prod_{G/H\in\pi_{0}\cO(G)} \mathbb P_{\mathcal E}\big(EW_G(H)\timesover{W_G(H)} X^H\big).
\]
\end{theorem}

\begin{proof}
Propositions~\ref{prop:FixedofFrees} and~\ref{prop:AutnGH} give the
stated splitting in the category of spaces.  We need to show that this
is a splitting of $\cE$-algebras; this is a consequence of
naturality.  

Recall from~\cite[6.25]{BHNinfty} that the cofree $\Ninfty$ operad
$F(EG,\cO)$ is a genuine $G$-$E_\infty$ operad for any $G$-operad
$\cO$.  Since the cofree functor is the right adjoint to the forgetful
functor to trivial operads, the unit gives us a natural map of
operads 
\[
\cO \to F(EG, i^*_e \cO)
\]
which provides an explicit model for the terminal map in the homotopy
category of $\Ninfty$ operads.  For convenience, we write $\cO^{gen}$
for $F(EG, i^*_e \cO)$ in what follows.

On the other hand, $\cO^{G}$ is an $E_\infty$
operad~\cite[B.1]{BHNinfty}, and so the inclusion gives rise to a map 
\[
\cO^{G} \to \cO
\]
of $\Ninfty$ operads where $\cO^{G}$ is given the trivial
$G$-action.

Putting this together, we see that the composite
\begin{equation}\label{eqn:OperadFixedPoints}
\big(\FreeO(X)\big)^G\to \big(\mathbb P_{\cO^{gen}}(X)\big)^G
\end{equation}
is a map of $\cO^{G}$-algebras.  The equivariant Barratt-Priddy-Quillen
theorem proved as~\cite[6.12]{GuillouMay} shows that we have an
equivalence of the desired form for $F(EG, i^*_e \cO)$.  In particular, the
projection maps  
\[
q_H\colon \big(\mathbb P_{\cO^{gen}}(X)\big)^G\to \mathbb P_{\mathcal E}(EW_G(H)\timesover{W_G(H)} X^H)
\]
are all maps of $\cO^{G}$-algebras.  Composing these maps
with the map in Equation~\ref{eqn:OperadFixedPoints} then shows that
the projection maps 
\[
q_H\colon \big(\FreeO(X)\big)^G\to \mathbb P_{\mathcal E}(EW_G(H)\timesover{W_G(H)} X^H)
\]
are all maps of $\cO^{G}$-algebras, and hence the splitting from
Proposition~\ref{prop:FixedofFrees} is one of $\cO^{G}$-algebras.
Naturality in maps of operads and the usual product trick then allows
us to convert this to a splitting of $\cE$-algebras, as desired.  
\end{proof}

We will now use the fact that for the free $\cO^{G}$-algebra functor
$\FreeO[\cO^{G}]$, the natural map  
\[
\FreeO[\cO^{G}] (X \vee Y) \to \FreeO[\cO^{G}] X \times \FreeO[\cO^{G}] Y
\]
is a weak equivalence~\cite[6.9]{BasterraMandell} and more generally
for cofibrant $\cO^{G}$-algebras, the natural map
\begin{equation}\label{eq:noneqstab}
X \vee_{\cO^{G}} Y \to X \times Y
\end{equation}
is a weak equivalence~\cite[6.8]{BasterraMandell}.

In particular, we can immediately deduce the following result
about naturality of the splitting in maps of operads.

\begin{corollary}\label{cor:ChangeofOperads}
If $\cO\to\cO'$ is a map of $\Ninfty$ operads, then for any $G$-CW complex $X$, the induced map
\[
\big(\FreeO(X)\big)^G\to \big(\mathbb P_{\cO'}(X)\big)^G
\]
is the inclusion of an $E_\infty$ direct summand.
\end{corollary}

As another corollary of the proof, we can maintain homotopical control
on the fixed points:

\begin{corollary}\label{cor:fixcof}
For an $\Ninfty$ operad $\cO$, the fixed points
\[
\big(\FreeO(X)\big)^G\simeq \prod_{G/H\in\pi_{0}\cO(G)} \mathbb P_{\mathcal E}\big(EW_G(H)\timesover{W_G(H)} X^H\big).
\]
have the homotopy type of a cofibrant $\cO^G$-algebra.
\end{corollary}

We can use these results to express the decomposition in terms of the
free algebra on the summands, as follows.  For $H\subset G$ such that
$G/H$ is admissible for $\cO$, let $\Gamma_{G/H}$ be the graph subgroup of $G\times\Sigma_{|G/H|}$ corresponding to an ordering of the points of $G/H$. 
Our description of the $G$-fixed points of $\FreeO(X)$ shows that
there is an inclusion 
\[
EW_{G}(H)\timesover{W_{G}(H)}X^{H}\xrightarrow{n_{G/H}} \FreeO(X)^{G}
\]
corresponding to the map on homotopy colimits induced by the inclusion
of the orbit $G\times\Sigma_{|G/H|}/\Gamma_{G/H}$ into the family
$\cF_{|G/H|}$.  By the free-forget adjunction, this gives us a map of
$\cO^{G}$-algebras 
\[
\mathbb P_{\cO^{G}}\big(EW_{G}(H)\timesover{W_{G}(H)}X^{H}\big)\xrightarrow{n_{G/H}} \FreeO(X)^{G}.
\]
The following proposition now follows immediately from
Theorem~\ref{thm:BPQ} and equation~\eqref{eq:noneqstab}.

\begin{proposition}\label{prop:InclusionofSummands}
The map
\[
\prod n_{G/H}\colon \mathbb P_{\cO^G}\!\left(\coprod_{G/H\in\pi_{0}\cO(G)}EW_{G}(H)\timesover{W_{G}(H)}X^{H}\right)\to \FreeO(X)^{G}
\]
is a weak equivalence of $\cE$-algebras.
\end{proposition}

As a corollary, we deduce a first weak kind of stability for the
category of $\Ninfty$-spaces.

\begin{corollary}\label{cor:coprodprod}
Let $X$ and $Y$ be cofibrant $\cO$-algebras in the standard model
structure.  Then the natural map 
\[
X \vee_{\cO} Y \to X \times Y 
\]
is an equivalence of $\cO$-algebras.
\end{corollary}

\begin{proof}
First, observe that using the bar construction equivalence
$B(\FreeO, \FreeO, X) \to X$ we can reduce as in the argument
for~\cite[6.8]{BasterraMandell} to the case where we are considering
the natural map 
\[
\FreeO X \vee_{\cO} \FreeO Y \to 
\FreeO X \times \FreeO.
\]
It suffices to consider the $G$-fixed points, and we can now apply
Theorem~\ref{thm:BPQ} and use the fact that $\FreeO X \vee_{\cO}
\FreeO Y \cong \FreeO(X \vee Y)$; we have that
\[
\big(\FreeO(X \vee Y)\big)^G\simeq \prod_{G/H\in\pi_{0}\cO(G)} \mathbb
P_{\mathcal E}\big(EW_G(H)\timesover{W_G(H)} (X \vee Y)^H\big).
\]
Since coproducts commute with fixed-points and orbits, we can write
the righthand side as 
\[
\prod_{G/H\in\pi_{0}\cO(G)} \mathbb P_{\mathcal
  E}\big(EW_G(H)\timesover{W_G(H)} X^H\big) \vee_{\cE} P_{\mathcal
  E}\big(EW_G(H)\timesover{W_G(H)} Y^H\big).
\]
Using the non-equivariant comparison of coproduct and product for
$\cE$-spaces (equation~\eqref{eq:noneqstab}), we can rewrite again as 
\[
\prod_{G/H\in\pi_{0}\cO(G)} \mathbb P_{\mathcal
  E}\big(EW_G(H)\timesover{W_G(H)} X^H\big) \times 
\prod_{G/H\in\pi_{0}\cO(G)}
P_{\mathcal
  E}\big(EW_G(H)\timesover{W_G(H)} Y^H\big),
\]
which by Theorem~\ref{thm:BPQ} is precisely
\[
\FreeO(X) \times \FreeO(Y).
\]
\end{proof}

\subsubsection{The based case}

We can deduce the corresponding splitting in the based case from the
unbased statement.  Recall that we have a natural isomorphism of monads
$\redFreeO((-)_+) \cong \FreeO(-)$.  Moreover, for any based $G$-space
$X$, there is a map $S^0 \to X_{+}$ specified by sending the non-basepoint
of $S^0$ to the basepoint of $X$.  If $X$ has a nondegenerate
basepoint (which we can assume without loss of generality), then there
is a cofiber sequence of based $G$-spaces $S^0 \to X_+ \to X$.  Since
$\redFreeO$ preserves cofiber sequences and so does the
righthand side of Theorem~\ref{thm:BPQ}, we can conclude the following
result in the based case.

\begin{theorem}\label{thm:BasedBPQ}
The norm maps $n_{G/H}$ for $G/H$ admissible induce a natural
weak equivalence of $\cO^{G}$-algebras 
\[
\big(\redFreeO(X)\big)^{G}\simeq \prod_{G/H\in\pi_{0}\cO(G)} \tilde{\mathbb P}_{\cO^{G}}(EW_{G}(H)_{+}\smashover{W_{G}(H)} X^{H}).
\]
\end{theorem}

Delooping this gives us the tom Dieck splitting.

\begin{proof}[Proof of Theorem~\ref{thm:OsptomDieck}]
The definition of $\Sp^{G}_{\cO}$ implies that we can compute the
fixed points of a suspension spectrum in terms of the zero-space.
Specifically, Using Lemma~\ref{lem:tensorfixed},
Corollary~\ref{cor:omega} implies that we can deduce the incomplete
tom Dieck splitting from Theorem~\ref{thm:BPQ}.
\end{proof}

\section{Equivariant stability}\label{sec:OSpectraProperties}
\subsection{The Wirthm\"uller Isomorphism in $\cO$-spectra}

The infinite loop space version of the Wirthm\"uller isomorphism is that for every subgroup $H\subset G$ and for every \(\cO^{gen}\)-algebra \(B\) in based $H$-spaces, we have a natural weak equivalence
\[
G_{+}\otimes_{H}B\xrightarrow{\simeq} \Map^{H}(G_{+},B).
\]

This is true much more generally. 

\begin{notation}
Let 
\[
\delta\colon i_{H}^{\ast}G_{+}\to H_{+}
\]
be the based $H$-map
\[
\delta_{eH}(g)=
\begin{cases}
g & g\in H \\
+ & g\notin H.
\end{cases}
\]
We will also use $\delta$ to denote the natural map of based $H$-spaces
\[
i_{H}^{\ast}(G_{+}\smashover{H}X)\to H_{+}\smashover{H}X\cong X.
\]
\end{notation}

\begin{theorem}\label{thm:RightAdjoint}
Let $H\subset G$ be such that $G/H$ is an admissible $\cO$-set. Then the natural map of $\cO$-algebras
\[
\redFreeO(G_{+}\smashover{H}X)\xrightarrow{\WirtO}\Map_{H}\big(G_{+},\redFreeO[\iota^*_H
  \cO](X)\big)
\]
adjoint to 
\[
i_{H}^{\ast}\redFreeO(G_{+}\smashover{H}X)\cong \redFreeO[\iota^*_H \cO](i_{H}^{\ast}G_{+}\smashover{H}X)\xrightarrow{\redFreeO(\delta)} \redFreeO(X)
\]
is a weak equivalence in the standard model structure.
\end{theorem}
The proof of Theorem~\ref{thm:RightAdjoint} is a somewhat elaborate
excursion through finite group theory, so we briefly postpone it to
state and prove our desired result: a Wirthm\"uller isomorphism for
the $\cO$-Spanier-Whitehead category. 

\begin{theorem}\label{thm:SWwirth}
For $H$ an admissible $\cO$-set, the functor
\[
X\mapsto G_{+}\smashover{H}X
\]
is both the left and right adjoint to the forgetful functor $i_{H}^{\ast}$ on $\SW^{H}_{\cO}$.
\end{theorem}
\begin{proof}
Since induction is the left adjoint in $G$-spaces, and since the mapping sets in $\SW^{G}_{\cO}$ can be computed as just mapping sets in $G$-spaces, induction is the left adjoint in $\SW^{H}_{\cO}$.

Since limits in $\cO$-algebras are formed in the underlying category, we know that $\Map_{H}(G_{+},-)$ is also the right adjoint to the forgetful functor 
\[
i_{H}^{\ast}\colon \cO\mhyphen\Alg^{G}\to\cO\mhyphen\Alg^{H}.
\]
Theorem~\ref{thm:RightAdjoint} gives us a natural equivariant weak equivalence for any $H$-space $Y$
\[
\Omega\redFreeO(G_{+}\smashover{H}\Sigma
Y)\to\Omega\Map_{H}\big(G_{+},\redFreeO[\iota^*_H \aO](\Sigma Y)\big)\cong 
\Map_{H}\big(G_{+},\Omega\redFreeO[\iota^*_H \aO](\Sigma Y)\big).
\]
Thus for any finite $G$-CW complex $X$, we have a natural
isomorphism
\begin{align*}
\{X,G_{+}\smashover{H}Y\}_{\cO}^{G}
&=\cO\mhyphen\Alg^{G}\Big(\Omega\redFreeO(\Sigma X),\Omega\redFreeO(G_{+}\smashover{H}\Sigma Y)\Big) \\
&\cong \cO\mhyphen\Alg^{G}\Big(\Omega\redFreeO(\Sigma
X),\Map_{H}\big(G_{+},\Omega\redFreeO[\iota^*_H \aO](\Sigma Y)\big)\Big)\\
&\cong \cO\mhyphen\Alg^{H}\Big(i_{H}^{\ast}\Omega\redFreeO(\Sigma
X),\Omega\redFreeO[\iota^*_H \aO](\Sigma Y)\Big)\\
&\cong \{i_{H}^{\ast}X,Y\}_{\cO}^{H},
\end{align*}
where the last isomorphism is the obvious one that the free algebra in $G$-spaces on $X$ restricts to the free algebra in $H$-spaces on $i_{H}^{\ast}X$.
\end{proof}

Finally, we can conclude the Wirthmuller isomorphism in $\Sp_G^{\cO}$;
our argument is analogous in spirit to the formal criterion given
in~\cite{FauskHuMay}.  For any subgroup $H \subset G$ and
$\iota^\ast_H\cO$-spectrum $X$, there is a natural map
\[
G_+ \otimes_H X \to F_H(G, X)
\] 
obtained as the adjoint of the natural map $X \to
\iota_H^\ast F_H(G,X)$ that takes $x$ to the map $\delta(-) x$.

\begin{theorem}\label{thm:OWirth}
Let $B$ be an $\iota_H^\ast \cO$-spectrum.  If \(G/H\) is an admissible
\(G\)-set, then there is a natural isomorphism in $\Ho(\Sp^{G}_{\cO})$  
\[
G_+ \otimes_{H} B \cong F_H(G, B).
\]
\end{theorem}

\begin{proof}
Since both sides commute with sifted homotopy colimits, it suffices to
show that the map is an isomorphism for the generating objects of
$\Sp^{H}_{\cO}$; by Proposition~\ref{prop:generators}, these are the
suspension spectra of $H/K$ where $K$ ranges over the subgroups of
\(H\). These are in particular suspension spectrum \(H\)-CW complexes,
and the result now follows from Theorem~\ref{thm:SWwirth}.
\end{proof}

\subsection{Proof of Theorem~\ref{thm:RightAdjoint}}
For self-containedness, we will prove in this subsection that the natural map
\begin{equation}\label{eqn:AdjointMap}
\redFreeO(G_{+}\smashover{H}X)\xrightarrow{\WirtO}\Map_{H}(G_{+},\redFreeO(X)\big)
\end{equation}
is a weak equivalence. If we restrict both sides to a proper subgroup $K\subset G$, then we have the map
\[
\redFreeO(i_{K}^{\ast}G_{+}\smashover{H}X)\xrightarrow{\WirtO} \Map_{H}\big(i_{K}^{\ast}G_{+},\redFreeO(X)\big),
\]
so by induction on the subgroup lattice, it will suffice to show that the map $\redFreeO(\delta)$ is an equivalence on $G$-fixed points. We will first analyze the fixed points, then use our usual comparison trick to deduce that the map is a weak equivalence.

The $G$-fixed points of the target is elementary, since we are mapping into something coinduced. 
\begin{proposition}\label{prop:FixedofTarget}
We have a natural weak equivalence of $\cO^{tr}$-algebras
\[
\Map_{H}(G_{+},\redFreeO(X)\big)^{G}\simeq \prod_{H/K\in\pi_{0}\cO(H)}\redFreeO\big(EW_{H}(K)_{+}\smashover{W_{H}(K)}X^{K}\big)
\]
\end{proposition}
\begin{proof}
For any $\cO$-space, the diagonal map provides a natural weak equivalence of $\cO^{tr}$-algebras
\[
\redFreeO(X)^{H}\xrightarrow{\simeq}\Map_{H}(G_{+},\redFreeO(X)\big)^{G}.
\]
The result then follows from Theorem~\ref{thm:BasedBPQ}.
\end{proof}

Theorem~\ref{thm:BasedBPQ} also describes the fixed points of the source as an $\cO^{tr}$-algebra.

\begin{proposition}\label{prop:FixedofSource}
We have a natural weak equivalence of $\cO^{tr}$-algebras
\[
\redFreeO(G_{+}\smashover{H}X)^{G}\simeq \prod_{G/K\in\pi_{0}\cO(G)}\redFreeO\big(EW_{G}(K)_{+}\smashover{W_{G}(K)}(G_{+}\smashover{H}X)^{K}\big).
\]
\end{proposition}

At this point, it is not even clear that the two products we have to
compare run over the same indexing set. We begin
simplifying. Additionally, we must pass from the Weyl groups in $G$ to
those in $H$. With a miraculous result about finite groups, all of this
is possible. 

We begin by analyzing the $K$-fixed points of $G_{+}\smashover{H}X$. First note that the double coset decomposition of $G$ gives an identification
\[
i_{K}^{\ast}\left(G_{+}\smashover{H}X\right)\cong \bigvee_{KgH\in K\backslash G/H} KgH_{+}\smashover{H} X,
\]
where we use $KgH$ both for the double coset as an equivalence class and as a $K$-$H$-biset. As a $K$-set, we have a natural (in $X$) $K$-equivariant homeomorphism 
\[
KgH_{+}\smashover{H}X\cong K_{+}\smashover{(K\cap gHg^{-1})} c_{g}^{\ast}X,
\]
where $c_{g}^{\ast}\colon\Top^{H}\to\Top^{gHg^{-1}}$ is the pull-back along the isomorphism $gHg^{-1}\cong H$. 

\begin{proposition}\label{prop:KFixedPoints}
We have a natural homeomorphism
\[
\big(G_{+}\smashover{H}X\big)^{K}\cong \bigvee _{\substack{KgH\in K\backslash G/H,\\ g^{-1}Kg\subset H}} X^{g^{-1}Kg}.
\]
\end{proposition}
We need this with the Weyl action of $K$ in $G$ to understand the homotopy orbits. Here it acts both on the indexing set and on the fixed points. 
%
In particular, we have an action on the double cosets.

\begin{proposition}
The normalizer of $K$ in $G$ acts on the set of double cosets by
\[
n\cdot Kg'H = nKg'H=Kng'H.
\]
This gives an obvious action of the Weyl group. 
\end{proposition}
\begin{proof}
The $H$ term is largely immaterial: since we are multiplying any double coset $Kg'H$ on the left by $n$ (either before or after the $K$), replacing $g'$ with $g'h$ for any $h\in H$ yields the same value. Similarly, since $g\in N_{G}(K)$, we know that if we replace $g'$ by $kg'$ for some $k\in K$, then
\[
Knkg'H=Kk'ng'H=Kng'H,
\]
where $k'=nkn^{-1}\in K$. This is therefore well defined. It is also obviously a left group action. Since $kK=K$ for all $k\in K$, we see that $K$ acts trivially, giving the final part.
\end{proof}

In general, the stabilizers could be a little ugly. However, we consider only those double cosets $KgH$ for which $g^{-1}Kg\subset H$.

\begin{definition}
Let $(K:G:H)$ denote the set of double cosets $KgH$ such that $g^{-1}Kg\subset H$.
\end{definition}

\begin{proposition}
The subset $(K:G:H)$ is an $W_{G}(K)$-equivariant subset of $K\backslash G/H$.
\end{proposition}
\begin{proof}
Let $w\in W_{G}(K)$ be represented by $n\in N_{G}(K)$, and let $KgH\in (K:G:H)$. Then since $n$ normalizes $K$,
\[
(ng)^{-1}Kng=g^{-1}n^{-1}Kng=g^{-1}Kg.
\]
\end{proof}

\begin{proposition}\label{prop:WeylAction}
Let $g$ be such that $g^{-1}Kg\subset H$. Then the stabilizer in $W_{G}(K)$ of $KgH$ is $W_{gHg^{-1}}(K)$ which under the conjugation isomorphism is $W_{H}(g^{-1}Kg)$.
\end{proposition}
\begin{proof}
Let $n\in N_{G}(K)$ stabilize $KgH$. Equivalently, there are $k\in K$ and $h\in H$ such that 
\[
ng=kgh.
\]
This means that $ghg^{-1}\in N_{G}(K)$, and as elements of the Weyl group, 
\[
[n]=[ghg^{-1}].
\]
The second part is obvious.
\end{proof}

\begin{corollary}\label{cor:WeylActionOnFixed}
We have an $W_{G}(K)$-equivariant homeomorphism
\[
\big(G_{+}\smashover{H}X\big)^{K}\cong
\bigvee_{[KgH]\in (K:G:H)/W_{G}(K)} W_{G}(K)_{+}\smashover{W_{gHg^{-1}}(K)}c_{g}^{\ast}X^{g^{-1}Kg}.
\]
\end{corollary}

Putting this together with Proposition~\ref{prop:KFixedPoints} then gives the following.

\begin{corollary}\label{cor:KFixedofSource}
We have a natural weak equivalence of $\cO^{tr}$-algebras
\[
\redFreeO(G_{+}\smashover{H}X)^{G}\simeq
\prod_{\substack{G/K\in\pi_{0}\cO(G),\\ K<H}} \prod_{[KgH]\in (K:G:H)/W_{G}(K)}\redFreeO\big(EW_{H}(g^{-1}Kg)_{+}\smashover{W_{H}(g^{-1}Kg)}X^{g^{-1}Kg}\big).
\]
Here $K<H$ denotes the subconjugacy relation.
\end{corollary}

This actually brings us much closer to finishing our proof. The indexing set for the first product is the $G$-conjugacy classes of subgroups $K$ such that $G/K$ is admissible and $K$ is subconjugate to $H$. Since $G/H$ is admissible, by assumption, this is the same as those $G$-conjugacy classes of subgroups $K$ such that $H/K$ is admissible.
\begin{proposition}
If $G/H$ is admissible and if $K\subset H$, then $G/K$ is admissible if and only if $H/g^{-1}Kg$ is admissible for all $g^{-1}Kg\subset H$.
\end{proposition}
\begin{proof}
Since $G/H$ is admissible, and since admissible sets are closed under self-induction, if $H/K$ is admissible, then 
\[
G/K\cong G\timesover{H}H/K
\]
is admissible. For the other direction, if $G/K$ is admissible, then $G/g^{-1}Kg$ is admissible for all $g\in G$, since these are isomorphic $G$-sets. In particular, for all $g$ such that $g^{-1}Kg\subset H$, $G/g^{-1}Kg$ is an admissible $G$-set. Restricting to $H$, this gives
\[
i_{H}^{\ast}G/g^{-1}Kg=T_{g}\amalg H/g^{-1}Kg
\]
for some finite $H$-set $T_{g}$. Since admissible sets are closed under finite limits, we conclude that $H/g^{-1}Kg$ is admissible.
\end{proof}

\begin{corollary}
The indexing set for the first product in Corollary~\ref{cor:KFixedofSource} is the set of admissible $H$-orbits modulo conjugation in $G$.
\end{corollary}

This explains how the two sides of the map in Equation~\ref{eqn:AdjointMap} for Theorem~\ref{thm:RightAdjoint} are packaging the data: we take the conjugacy classes of subgroups of $H$ and group them according to conjugacy in $G$.

\begin{definition}
Let $H, K\subset G$. Let
\[
(K)_{G;H}=\big\{g^{-1}Kg \mid g\in G, g^{-1}Kg \subset H\big\}/H\text{-conjugacy}
\]
be the set of $H$-conjugacy classes of $G$-conjugates of $K$ which sit in $H$.
\end{definition}
In other words, $(K)_{G;H}$ is what we get when we take the conjugacy class of $K$, throw out all the terms which are not subgroups of $H$ and then work up to $H$-conjugacy.

\begin{proposition}
The map $\theta\colon (K:G:H)\to (K)_{G;H}$ defined by
\[
\theta(KgH)=g^{-1}Kg
\]
is Weyl equivariant (where the target has a trivial Weyl action) and induces a bijection
\[
(K:G:H)/W_{G}(K)\cong (K)_{G;H}.
\]
\end{proposition}
\begin{proof}
The map $\theta$ is obviously onto if it is well-defined. However, since we are considering only $H$-conjugacy classes of subgroups of $H$, replacing $g$ by $kgh$ results in an $H$-conjugate subgroup, and the map is well-defined.

For the second part, if $n\in N_{G}(K)$, then as before
\[
(ng)^{-1}Kng=g^{-1}Kg,
\]
and hence $\theta$ factors through the Weyl orbits. On the other hand, if 
\[
\theta(KgH)=\theta(Kg'H),
\]
then there is an $h\in H$ such that 
\[
g^{-1}Kg=h^{-1}{g'}^{-1}Kg'h=(g'h)^{-1}K(g'h).
\]
Rearranging, this means that $n=g'hg^{-1}\in N_{G}(K)$, and by construction, 
\[
n\cdot KgH=K(g'hg^{-1})gH=Kg'hH=Kg'H.
\]
\end{proof}

\begin{corollary}
The left- and right-hand sides of the $G$-fixed points of Equation~\ref{eqn:AdjointMap} are abstractly equivalent as $\cO^{tr}$-algebras.
\end{corollary}

We need slightly more, in that we need that the map $\WirtO^{G}$ from
equation~\eqref{eqn:AdjointMap} is an equivalence of $\cO^{tr}$-algebras. Here we again compare to a larger operad in which we know the answer. Recall from Corollary~\ref{cor:ChangeofOperads} that if $\cO'$ is any larger $\Ninfty$ operad, then the fixed points of the free $\cO$-algebra sit as the obvious $\cO^{tr}$-summand of the fixed points of the free $\cO'$-algebra. The map in Equation~\ref{eqn:AdjointMap} fits into a commutative square
\begin{center}
\begin{tikzcd}
{\redFreeO(G_{+}\smashover{H}X)}
        \arrow[r, "\WirtO"]
        \ar[d]
        &
{\Map_{H}\big(G_{+},\redFreeO(X)\big)}
        \ar[d]
        \\
{\tilde{\mathbb P}_{\cO^{'}}(G_{+}\smashover{H}X)}
                \ar[r, "\omega_{\cO',H}"]
                &
{\Map_{H}\big(G_{+},\tilde{\mathbb P}_{\cO'}(X)\big).}
\end{tikzcd}
\end{center}
Taking fixed points of the bottom row and using that the vertical maps are the inclusion of summands, we see that the map in Equation~\ref{eqn:AdjointMap} for $\cO'$ can be written as
\begin{equation}\label{eqn:AdjointSummands}
\big(\redFreeO(G_{+}\smashover{H}X)\big)^{G}\vee_{\cO^{tr}} E_{s}\xrightarrow{\left[\begin{matrix}
\WirtO^{G} & ? \\
0 & ?
\end{matrix}\right]} 
\Big({\Map_{H}\big(G_{+},\redFreeO(X)\big)}\Big)^{G}\vee_{\cO^{tr}} E_{t},
\end{equation}
where $E_{s}$ and $E_{t}$ are the complementary summands. Here
$\omega_{\cO',H}^{G}$ is a ``square matrix'' in the sense that the source and
target are abstractly isomorphic $\cO^{tr}$-algebras. 

Now consider any universe $U$ in which $G/H$ equivariantly embeds and
for which $\ccD(U)$ is bigger than $\cO$ (at worst, one can take a
complete universe). Then Adams' original argument for the Wirthmuller
isomorphism~\cite[\S 5]{Adams} shows that the map 
\[
\tilde{\mathbb
  P}_{\ccD(U)}(G_{+}\smashover{H}X)\to\Map_{H}\big(G_{+},\tilde{\mathbb
  P}_{\ccD(U)}(X)\big) 
\]
is an equivariant weak equivalence. In particular, it is a weak
equivalence on fixed points. The matrix form of this map in
Equation~\ref{eqn:AdjointSummands} then shows that $\WirtO^{G}$ is a
weak equivalence, as desired.

\section{Comparisons}

The purpose of this section is to compare the category of
$\cO$-spectra for certain $\Ninfty$ $G$-operads $\cO$ to the usual
categories of equivariant orthogonal $G$-spectra.  Specifically, we
explain the proof of the following theorem.

\begin{theorem}\label{thm:compare}
Let $\cO$ be an $\Ninfty$ operad which is equivalent to the
equivariant Steiner operad for some universe $U$.  Then there is an
equivalence of $\infty$-categories between $N(\Sp^{G}_{\cO})[W^{-1}]$
and $N(\Sp^{G})[W_{U}^{-1}]$; i.e., between the $\cO$-stable category
and the category of orthogonal $G$-spectra with the stable
equivalences determined by $U$.
\end{theorem}

Here by the equivariant Steiner operad for a universe $U$, we mean the
following.  First, recall the definition of the Steiner operad $\cK_V$
for a finite-dimensional real inner product $G$-space $V$.  A Steiner
path is a map $h \colon I \to \Map(V,V)$ such that $I$ lands in the
distance-reducing linear embeddings and $h(1)$ is the identity map.
The $n$th space of the Steiner operad consists of tuples $(h_1,
\ldots, h_n)$ of Steiner paths such that the embeddings $\{h_i(0)\}$
have disjoint image.  The symmetric group acts by permuting the
tuples, $G$ acts by conjugation on embeddings, and the composition is
determined by the pointwise composition.  The homotopy type of $\cK_V$
can be described in terms of an equivariant configuration space.

For $V \subset W$, there is an evident inclusion $\cK_V(n) \to
\cK_W(n)$ that is equivariant for the $G \times \Sigma_n$ action.
Moreover, these inclusions assemble into a map of operads $\cK_V \to
\cK_W$.  Therefore, for a universe $U$, we define the equivariant
Steiner operad for $U$ as
\[
\cK_U = \colim_{V \in U} \cK_V,
\]
the colimit over the collection of finite-dimensional subspaces of $U$
with morphisms inclusions.

The proof of Theorem~\ref{thm:compare} depends on the following
result, which is essentially folklore but has not previously appeared
in the literature.

\begin{theorem}\label{thm:opcompare}
There is a zig-zag of connective Quillen equivalences between the
category $\GTop[\bO]$ with the grouplike model structure and the
category $\Sp_G$ of orthogonal $G$-spectra with the stable equivalences
determined by $U$.
\end{theorem}

Here recall from the discussion
preceding~\cite[0.10]{MandellMaySchwedeShipley} that a connective
Quillen equivalence is a Quillen equivalence which restricts to an
equivalence on the full subcategories of the respective homotopy
categories spanned by the connective objects.

We first explain the proof of Theorem~\ref{thm:compare} given
Theorem~\ref{thm:opcompare}.  This is an immediate consequence of the
following lemma.  Here let $\Sp_G^{c}$ denote the subcategory of
connective objects in $\Sp_G$ and denote by $\Sp(\aC)$ the
$\infty$-category of spectrum objects in an $\infty$-category
$\aC$~\cite[\S 1.4.2]{LurieHA}.

\begin{lemma}
There is an equivalence of $\infty$-categories
\[
\Sp(N(\Sp^c_G)[W_{U}^{-1}]) \htp N(\Sp_G)[W_U^{-1}].
\]
\end{lemma}

\begin{proof}
The inclusion 
\[
N(\Sp_G^{c})[W_{U}^{-1}] \to N(\Sp_G)[W_U^{-1}]
\]
is exact and therefore induces a functor 
\[
\Sp(N(\Sp_G^{c})[W_{U}^{-1}]) \to \Sp(N(\Sp_G)[W_U^{-1}]).
\]
Computation of homotopy groups shows that this is an equivalence, and
since the canonical suspension spectrum functor
\[
N(\Sp_G)[W_U^{-1}] \to \Sp(N(\Sp_G)[W_U^{-1}])
\]
is an equivalence because $N(\Sp_G)[W_U^{-1}]$ is stable, the result
follows.
\end{proof}

We now discuss the proof of Theorem~\ref{thm:opcompare}.  Since $\cO$
is equivalent to the equivariant Steiner operad on $U$, we can
approach this using the basic strategy of the (equivariant analogue)
of~\cite[5.43]{ABGHR}.  Here we replace the specific discussion of the
operadic bar construction with the work of~\cite[\S
  6]{MerlingMayOsorno}.  The outline for this kind of result is given
in~\cite[\S 9]{MayRant}, and we now explain the verification of the
specific facts we need.  Recall that $\bQ$ denotes the monad
associated to the construction $\colim_V \Omega^V \Sigma^V X$ as $V$
varies over the indexing spaces in the universe.

\begin{enumerate}

\item The map of monads from $\bK_U \to \bQ$ is an equivariant group
  completion.  This follows from the work of
  Caruso-Waner~\cite{carusowaner}, as reviewed and stated
  in~\cite[1.11]{GuillouMay}.

\item The geometric realization of proper simplicial $G$-spaces and
  orthogonal $G$-spectra preserves levelwise weak equivalences.  This
  is standard (e.g., see~\cite[\S 4]{goodthh} for a discussion).

\item Finally, we need to know that $\Omega^V |X_\bullet| \cong
  |\Omega^V X_\bullet|$ for suitably connected $X$; as reviewed in the
  proof of~\cite[1.14]{GuillouMay}, this is due to Hauschild and a
  published source is~\cite{costenoblewaner}.

\end{enumerate}

It is natural to wonder about the analogues of the comparison results
in this section for arbitrary indexing systems.  As we discussed
briefly in the introduction, the subtleties of the situation
are illuminated by consideration of the definition of very special
$\Gamma$-$G$-spaces.  Let $T$ be a based finite $G$-set, and $X$ a
$\Gamma$-$G$-space; recall this is a functor from the category of
based finite sets to $G$-spaces.  Then $X$ is special if the natural
map 
\[
X(T) \to \Map(T, X(1^+))
\]
is an equivalence of $G$-spaces, where $X(T)$ denotes the value of $X$
on the underlying set of $T$ and $1^+$ denotes the one-point set with
a disjoint basepoint added.  We can of course ask for the special
condition to hold only for admissible based finite $G$-sets $T$; this
corresponds to indexing systems which we denote ``disk-like'', where
the relevant data is determined by what happens at $G/G$.  All little
disks operads give rise to disk-like indexing systems, but in fact
Rubin~\cite{Rubin2} shows that there exist disk-like indexing systems which
do not come from any little disks operad.  

There is a natural forgetful functor from $G$-$\Gamma$-spaces to
$H$-$\Gamma$-spaces, and one can use this to define versions of
special $\Gamma$-spaces in terms of conditions that hold for $H$-sets
for proper subgroups of $G$.  Generalizing this kind of condition, it
is possible to use excisive functors, spectral presheaves, and
equivariant symmetric spectra to model the $\cO$-stable category.  We
intend to return to these matters in future work.

\bibliographystyle{plain}

\bibliography{OSpectra}

\end{document}